\documentclass[final,3p,times]{elsarticle}

\usepackage{amssymb,amsmath,amsfonts,mathtools}
\usepackage{amsthm}
	\theoremstyle{remark}

	\theoremstyle{plain}

	\DeclareMathAlphabet{\mathbcal}{OMS}{cmsy}{b}{n}
	
\usepackage{algorithm}
\usepackage{algorithmicx,algpseudocode}
\usepackage{cases}

\usepackage[nodots]{numcompress}

\usepackage{hyperref}
\usepackage[]{setspace}

\biboptions{sort&compress}
\usepackage[normalem]{ulem}

\journal{arXiv}

\renewcommand{\vec}[1]{\mathbf{#1}}
\newcommand{\tensorTwo}[1]{\boldsymbol#1}
\newcommand{\tensorFour}[1]{\mathbb{#1}}

\newcommand{\blkMat}[1]{\boldsymbol{\mathsf{#1}}}
\newcommand{\blkVec}[1]{\boldsymbol{\mathsf{#1}}}

\newcommand{\Mat}[1]{\mathsf{#1}}
\renewcommand{\Vec}[1]{\mathsf{#1}}

\newcommand{\sub}[1]{_\mathit{#1}}

\usepackage{booktabs}
\usepackage{color, soul}
\usepackage[font=small]{caption}
\usepackage{subcaption}
\usepackage{arydshln}
\usepackage{upgreek}

\begin{document}

\begin{frontmatter}



\title{A macroelement stabilization for multiphase poromechanics}


\author[label1]{Julia T. Camargo\corref{cor1}}
\ead{jcamargo@stanford.edu}

\author[label2]{Joshua A. White}
\ead{jawhite@llnl.gov}

\author[label1]{Ronaldo I. Borja}
\ead{borja@stanford.edu}

\address[label1]{Department of Civil and Environmental Engineering,
                 Stanford University,
                 United States}

\address[label2]{Atmospheric, Earth and Energy Division,
                 Lawrence Livermore National Laboratory,
                 United States}
                 
\cortext[cor1]{Corresponding author}

\begin{abstract}

Strong coupling between geomechanical deformation and multiphase fluid flow appears in a variety of geoscience applications.  A common discretization strategy for these problems is a continuous Galerkin finite element scheme for the momentum balance equations and a finite volume scheme for the mass balance equations.   When applied within a fully-implicit solution strategy, however, this  discretization is not intrinsically stable.  In the limit of small time steps or low permeabilities, spurious oscillations in the pressure field, i.e. checkerboarding, may be observed.  Further, eigenvalues associated with the spurious modes will control the conditioning of the matrices and can dramatically degrade the convergence rate of iterative linear solvers.  Here, we propose a stabilization technique in which the balance of mass equations are supplemented with stabilizing flux terms on a macroelement basis.  The additional stabilization terms are dependent on a stabilization parameter.  We identify an optimal value for this parameter using an analysis of the eigenvalue distribution of the macroelement Schur complement matrix. The resulting method is simple to implement and preserves the underlying sparsity pattern of the original discretization.  Another appealing feature of the method is that mass is exactly conserved on macroelements, despite the addition of artificial fluxes.  The efficacy of the proposed technique is demonstrated with several numerical examples.
\end{abstract}

\begin{keyword}
poroelasticity \sep
reservoir simulation \sep
inf-sup stability  \sep
finite element method \sep
finite volume method



\end{keyword}

\end{frontmatter}


\allowdisplaybreaks


\section{Introduction}
\label{sec:intro}

In a variety of applications, it is useful to model the hydromechanical behavior of porous media infiltrated by one or more fluids---e.g. in geotechnical engineering \cite{teatini2006,borja2010,borja2012multi,camargo2016,cao2018,FaveroNeto2018,ghaffaripour2019,gholamiKorzani2018,mikaeili2018,navas2018,fusao2019,prassetyo2018,song2018,wang2018,WangWei2019,Yan2019a,Zhou2019}, hydrocarbon recovery \cite{zoback2010,Feng2019,Mashhadian2018,Peng2019,Semnani2016,Shiozawa2019,Yan2019b,Zhang2019,ZhaoYang2018,ZhaoYang2019}, and geologic carbon storage \cite{verdon2013, white2014,Jin2019}. Precise models should account for the tight interaction between solid deformation and fluid flow. The conceptual framework for modeling this coupled behavior is well established \cite{boer2000,coussy2005}, with Biot's work \cite{biot1941} providing a sound theoretical foundation. Computational methods for poromechanics, however, still pose many interesting challenges. In particular, this work focuses on numerical instabilities that may arise due to the discretization spaces chosen for the coupled fields.  

As a representative problem, we consider a model in which an elastic solid skeleton is saturated with two immiscible fluids.   We present a multiphase formulation for its relevance in many geoscience applications, though the single-phase formulation is a straightforward sub-case.  The behavior of the porous system is governed by a momentum balance equation for the mixture and mass balance equations for each of the fluids.   A fully-implicit time integration strategy is adopted, where all unknown fields are updated simultaneously in a monolithic manner \cite{white2016,white2018}.  A variety of finite-element and finite-volume based discretization strategies may be applied to these equations, each with certain advantages  \cite{settari1998,kim2011,prevost2014,castelletto2015,nordbotten2014,nordbotten2016,honorio2018,choo2018,sokolova2019}.  Of specific interest here, however, is a frequent choice: continuous trilinear interpolation for the displacement unknowns and element-wise constant fields for the pressure and saturation unknowns.   Such an interpolation results, for example, when applying a continuous Galerkin finite element scheme to the momentum balance equation, and a finite volume scheme to the mass balance equations \cite{settari1998,kim2011,prevost2014, white2018}.

This discretization works well in a variety of practical cases. The chosen interpolation spaces, however, can be problematic. In the limit of small time steps or low permeabilities, undrained deformation can occur. The fluid mass balance equations impose an incompressibility constraint on the deformation field.  Like many other constrained problems---e.g. Stokes flow or incompressible elasticity---these divergence constraints can create numerical instabilities if the discrete approximations for the field variables do not satisfy the Ladyzhenskaya-Babu\v{s}ka-Brezzi (LBB or inf-sup) condition \cite{babuska1973,brezzi1974}.   Unfortunately, the combination explored here is not intrinsically LBB-stable.  As a result spurious oscillations, i.e. checkerboarding, may be observed in the pressure field.  A less obvious, but equally important, symptom is a degradation in the convergence rate of iterative linear solvers.  Near-zero eigenvalues associated with spurious modes will control the conditioning of the system matrices, leading to poorly conditioned systems and increased iteration counts.  This latter issue can persist in regimes where checkerboarding is not visually apparent, and may thus go unnoticed by practitioners.

One way of circumventing these instabilities is to introduce a carefully-designed perturbation to the constraint equations.  The goal is to remove instabilities while maintaining an accurate approximation of the underlying PDEs. This is the basic rationale behind many stabilization techniques, including the Brezzi and Pitk\"aranta scheme \cite{brezzi1984}, the Galerkin Least-Squares approach \cite{hughes1986}, and the Polynomial Pressure Projection technique \cite{dohrmann2004}.  Various stabilization schemes have been proposed that devote particular attention to constant pressure elements \cite{fortin1990,bochev2006,Burman2007,hughes1987,silvester1990}, starting with early work in \cite{pitkaranta1985}.  In \cite{hughes1987} the idea of penalizing the pressure jump across inter-elements boundaries was introduced. An important modification of this method, the Local Pressure Jump (LPJ) stabilization, was developed in \cite{silvester1990} based on the macroelement concept.

The schemes above were primarily developed for fluid mechanics problems. Since then, many of these stabilization schemes have been successfully applied to poromechanics with single-phase flow \cite{white2008,preisig2011,choo2015,sun2013,berger2015,rodrigo2018,honorio2018}.   However, the study of stabilization procedures addressing multiphase problems is still incipient, with just a few studies available \cite{wan2002,truty2006}.  

This paper proposes a new stabilization technique in which the balance of mass equations are supplemented with stabilizing flux terms.  The resulting technique mimics the LPJ stabilization \cite{silvester1990,silvester1994,elman2005} in its basic design, but with suitable extensions to handle the poromechanical system of interest here.  The additional stabilization terms are dependent on a stabilization parameter that must be well chosen to suppress instabilities while not compromising solution accuracy.  We identify an optimal value for this parameter using an analysis of the eigenvalue distribution of the macroelement Schur complement matrix. The resulting method is simple to implement and preserves the underlying sparsity pattern of the original discretization.  Another appealing feature of the method is that mass is exactly conserved on macroelements.  
 

The paper is organized as follows. The governing equations and discretization scheme are introduced in Section~\ref{sec:model} and~\ref{sec:num_model}.   In Section~\ref{sec:incomp} we examine the behavior of this model in the undrained limit, in order to identify the source of spurious modes.  To fix this deficiency, our stabilization scheme is detailed in Section~\ref{sec:stabilized}. The resulting approach both treats spurious pore pressure oscillations and improves the conditioning of the system matrices. This is demonstrated through numerical examples presented in Section~\ref{sec:numerical_results}. Finally, concluding remarks are given in Section~\ref{sec:conclusions}.

\color{black}
\section{Governing Equations}
\label{sec:model}

We consider a multiphase poroelastic problem in which two immiscible fluids fill the voids of the porous, deformable solid skeleton. We focus on a displacement-saturation-pressure formulation, ignoring dynamic and non-isothermal effects. We further neglect capillary forces, meaning,  the wetting and non-wetting fluid phases have equal pressure inside the pore. This simplification is common in many reservoir-scale simulations \cite{Ertekin2000}.  The methods described here may be readily extended to include more sophisticated formulations.

The porous medium occupies a domain $\Omega \in \mathbb{R}^3$ over time interval $\mathcal{I} = (0,T]$. 	The unknown fields are the displacement of the solid $\vec{u} : \Omega \times \mathcal{I} \rightarrow \mathbb{R}^3$, the wetting fluid saturation $s : \Omega \times \mathcal{I} \rightarrow \mathbb{R}$, and the fluid pressure  $p : \Omega \times \mathcal{I} \rightarrow \mathbb{R}$. The initial/boundary value problem is governed by a linear momentum balance for the mixture and two mass balance equations for the wetting ($w)$ and non-wetting ($o$) fluids, respectively:
\begin{subequations}\label{eq:IBVP_global}
	\begin{align}
	\nabla \cdot \tensorTwo{\sigma'} -b \nabla p + \rho\vec{g} &= \vec{0} &\qquad \mbox{ in } \Omega \times \mathcal{I}, \label{eq:momentumBalanceS}\\
	%
	\dot{m}_{w} + \nabla \cdot \vec{w}_{w} - q_{w} &= 0 &\qquad \mbox{ in } \Omega \times \mathcal{I}, \label{eq:massBalanceW_S}	 \\
	%
	\dot{m}_{o}+ \nabla \cdot \vec{w}_{o} - q_{o} &= 0 & \qquad \mbox{ in } \Omega \times \mathcal{I}. \label{eq:massBalanceNW_S}
	\end{align}
\end{subequations}
In Eq.~\eqref{eq:momentumBalanceS}, the effective Cauchy stress depends on the symmetric gradient of the displacement field as
\begin{equation}
\tensorTwo{\sigma'} = \tensorFour{C} : \nabla^s \vec{u} ,
\end{equation}
where $\tensorFour{C}$ is the drained elasticity tensor.  Biot's coefficient $b=1-K/K_s$ may be calculate from $K$, the drained skeleton bulk modulus, and $K_s$, the intrinsic bulk modulus of the solid phase.  The mixture density $\rho$ is related to individual the phase densities---denoted by $\rho_\ell$ for $\ell=\{w,o\}$ and $\rho_s$---via the relationship 
\begin{equation}
\rho =  (1-\phi) \rho_s + \phi \rho_w s + \phi \rho_{o} (1-s) 
\end{equation}
where $\phi$ is the porosity.  Porosity changes are related to solid deformation and fluid pressure changes as
\begin{equation}
\dot \phi = b \nabla \cdot \dot{\tensorTwo{u}}+\frac{(b-\phi_0)}{K_s} \dot{p},
\end{equation}
which introduces a coupling between the momentum and mass balance equations.  Each fluid phase requires a density model $\rho_\ell(p)$, such as the simple linear model
\begin{equation}
\rho_\ell = \rho^0_\ell \left[1 + \frac{1}{K_\ell} (p-p^0) \right]
\end{equation}
with phase bulk modulus $K_\ell$ and reference density $\rho^0_\ell$ at a reference pressure $p^0$.   

In Eq.~\eqref{eq:massBalanceW_S} and Eq.~\eqref{eq:massBalanceNW_S}, $m_\ell$ is the mass per unit volume for the fluid phase $\ell=\{w,o\}$, with
\begin{subequations}
\begin{align}
 m_w &= \phi \rho_w s,\\
 m_{o} &= \phi \rho_{o} (1-s). 
 \end{align}
 \end{subequations}
The source terms $q_\ell$ are used to model well sources for injection and production of fluids, using a Peaceman well model \cite{Pea78,CheHuaMa06}.  The mass flux $\vec{w}_\ell= (\rho_\ell \vec{v}_\ell)$ is linked to the pore pressure field via the generalized Darcy's law as
\begin{equation}\label{eq:darcy}
\vec{v}_\ell = - \lambda_\ell \tensorTwo{\kappa} \cdot \nabla ( p + \rho_\ell g z ).
\end{equation}
The constitutive relation in \eqref{eq:darcy} defines the volumetric flux $\vec{v}_\ell$ using the phase mobility $\lambda_\ell = k_{r\ell}/\mu_\ell$, the viscosity $\mu_\ell$, and the relative permeability $k_{r\ell}$.  Specific relationships for viscosity $\mu_\ell = \mu_\ell(p)$ and relative permeability $k_{r\ell} = k_{r\ell} (s)$ must be defined for the fluids and porous medium under consideration.  Additionally, $\tensorTwo{\kappa}$ represents the absolute permeability tensor, $g$ the gravitational acceleration, and $z$ the elevation above a datum.

The domain boundary $\Gamma$ is decomposed into regions where Dirichlet and Neumann boundary conditions are specified, denoted by $\Gamma = \overline{\Gamma_{\vec{u}}^D \cup \Gamma_{\vec{u}}^N}$ for the momentum balance and  $\Gamma= \overline{\Gamma_f^D \cup \Gamma_f^N}$ for the mass balances.  	These divisions follow the overlap restriction $\Gamma_{\vec{u}}^D \cap \Gamma_{\vec{u}}^N = \Gamma_f^D \cap \Gamma_f^N = \emptyset$.  Specifically,
%
\begin{subequations}\label{eq:IBVP_boundary}
	\begin{align}
	&\vec{u} = \vec{0} & &\mbox{ on } \Gamma_{\vec{u}}^D \times \mathcal{I}, \label{eq:momentumBalanceS_DIR}\\
	&\tensorTwo{\sigma} \cdot \vec{n} = \bar{\vec{t}} & &\mbox{ on } \Gamma_{\vec{u}}^N \times \mathcal{I}, \label{eq:momentumBalanceS_NEU}\\
	&p =\bar{p} & &\mbox{ on } \Gamma_f^D \times \mathcal{I}, \label{eq:massBalanceS_DIR}\\    
	&s =\bar{s} & &\mbox{ on } \Gamma_f^D \times \mathcal{I}, \label{eq:massBalanceS_DIR_2}\\        
	&\vec{w}_{w} \cdot \vec{n} = 0 & &\mbox{ on } \Gamma_f^N \times \mathcal{I}, \label{eq:massBalanceW_S_NEU}\\           
	&\vec{w}_{o} \cdot \vec{n} = 0 & &\mbox{ on } \Gamma_f^N \times \mathcal{I},\label{eq:massBalanceNW_S_NEU}
	\end{align}   
	\end{subequations}
where the boundary conditions prescribing displacement \eqref{eq:momentumBalanceS_DIR}, total traction \eqref{eq:momentumBalanceS_NEU}, pore pressure \eqref{eq:massBalanceS_DIR}, wetting phase saturation \eqref{eq:massBalanceS_DIR_2}, wetting phase mass flux \eqref{eq:massBalanceW_S_NEU} and non-wetting phase mass flux \eqref{eq:massBalanceNW_S_NEU} are given. Here, $\vec{n}$ denotes the outer normal vector. Homogeneous conditions on the displacement and external fluxes were chosen here to simplify some notations below, but these can be easily relaxed.  

Initial conditions are specified as
	\begin{subequations}\label{eq:IBVP_initial}
	\begin{align}
	&\vec{u}(\vec{x}, 0) = \vec{u}_0 (\vec{x}) & &\mbox{ } (\vec{x}, t) \in (\Omega \times t=0), \label{eq:massBalanceDisp_IC}\\
	&s(\vec{x}, 0) = s_0 (\vec{x}) & &\mbox{ } (\vec{x}, t) \in (\Omega \times t=0), \label{eq:massBalanceSatW_IC}\\
	&p(\vec{x}, 0) = p_0 (\vec{x}) & &\mbox{ } (\vec{x}, t) \in (\Omega \times t=0). \label{eq:massBalancePres_IC}
	\end{align}
\end{subequations}
Note that the single-phase poromechanics model arises as a subcase of these general equations if one fixes either $s(\vec{x},t)=0$ or $s(\vec{x},t)=1$.  In this case, only one mass balance equation is required.

Clearly a number of modeling and constitutive assumptions are embedded in the formulation described above, but it remains a useful approximation for many subsurface applications.  This formulation also contains many of the salient mathematical features that may be encountered in other models used in practice.

\section{Discrete Formulation}
\label{sec:num_model}

\begin{figure}
\centering
\includegraphics[width=0.3\textwidth]{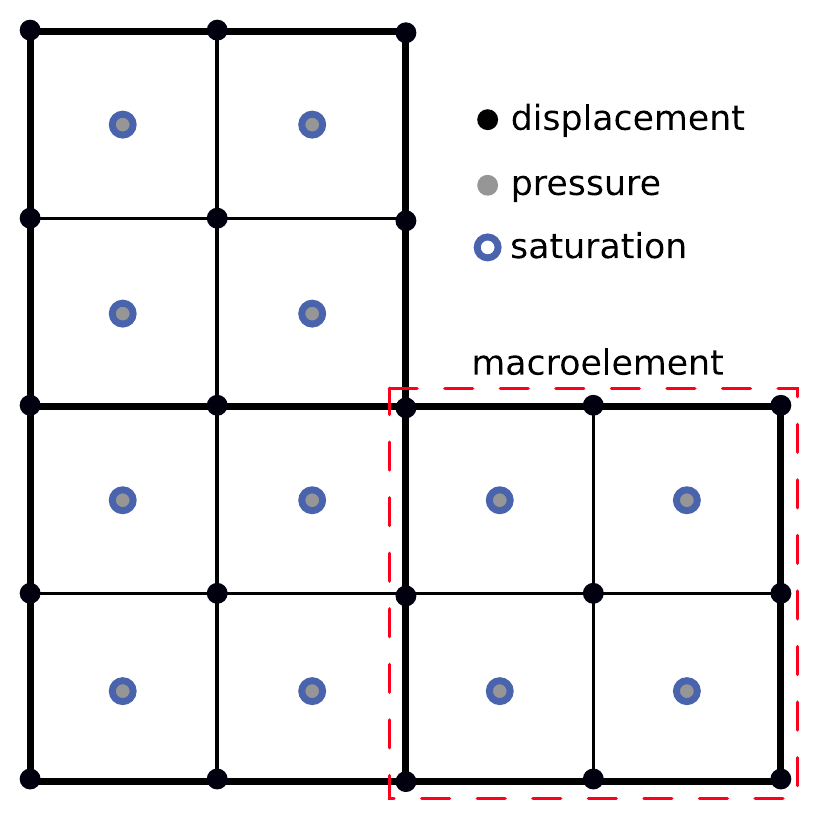}
\caption{Example mesh with nodal and cell-centered degrees-of-freedom.  Each element is assigned to a parent macroelement. }
\label{fig:geom}
\end{figure}

Figure~\ref{fig:geom} illustrates the basic geometry under consideration.  The domain $\Omega$ is partitioned into a computational mesh $\mathcal{T}^h$ made of non-overlapping elements $\{K_i\}$ such that $\Omega = \bigcup_{i=1}^{n_e} K_i$.  Every element face $f$ is assigned a unique unit normal vector $\vec{n}_f$. For our stabilization procedure, we further assume that these elements may be grouped into macroelements $\{M_i\}$ consisting of patches of eight hexahedra in 3D or four quadrilaterals in 2D.  This configuration is readily achieved by beginning with a coarse version of the mesh and applying one level of structured refinement. The discretization of the governing equations \eqref{eq:IBVP_global} is then obtained using a combined finite-element/finite-volume approach.   An extensive description of this formulation is presented in \cite{white2018}.  Here, we briefly summarize the key components, but refer the interested reader there for a more complete exposition.

Time integration relies on a fully implicit backward Euler scheme, with the time interval $\mathcal{I}$ divided into $n_{\Delta t}$ subintervals of length $\Delta t = (t_n - t_{n-1})$. We will use the notation $\Delta x=(x_n-x_{n-1})$ for other time-differenced quantities as well.  For the spatial discretization, we introduce the specific spaces
\begin{subequations}
\begin{align}
& \boldsymbol{\mathcal{Q}} := \left\lbrace \tensorTwo{v} \left|\right. \tensorTwo{v} \in [C^0(\Omega)]^3, \tensorTwo{v} ={\vec{0}} \text{ on } \Gamma_{\vec{u}}^{D},  \right.
 \left.\tensorTwo{v}_{\left| K \right.} \in [\mathbb{Q}_1(K)]^3 \; \forall K \in \mathcal{T}^h \right\rbrace, \label{eq:u_discrete_spaces} \\
& \mathcal{P} := \left\lbrace q \left|\right. q \in L^2(\Omega),  q_{\left| K \right.} \in \mathbb{P}_0(K) \; \forall K \in \mathcal{T}^h  \right\rbrace ,	\label{eq:s_discrete_spaces}
\end{align}
\end{subequations}
where $C^0(\Omega)$ and $L^2(\Omega)$ are the space of continuous and square Lebesgue-integrable functions on $\Omega$, respectively. $\mathbb{Q}_1(K)$ denotes the space of $d$-linear functions (trilinear in 3D or bilinear in 2D) and $\mathbb{P}_1(K)$ the space of constant functions on a given element $K$.

The discrete weak form of \eqref{eq:IBVP_global} reads: Find $\{ \vec{u}_n, s_n, p_n \} \in \boldsymbol{\mathcal{Q}} \times \mathcal{P} \times \mathcal{P}$ such that for time step $n = \{1, \ldots, n_{n_{\Delta t}} \}$
\begin{subequations} \label{eq:weak_form_discr}
	\begin{align}
	R_{u} =&
	\int_{\Omega} \nabla^{s} \tensorTwo{v} : \tensorTwo{\sigma}^\prime_n \, \mathrm{d} \Omega -
	\int_\Omega  \nabla \cdot \tensorTwo{v} b p_n \, \mathrm{d}\Omega 
	 - \int_\Omega \tensorTwo{v} \cdot \rho_n \vec{g} \, \mathrm{d}\Omega  
	 - \int_{\Gamma_{\vec{u}}^{N}} \tensorTwo{v} \cdot \bar{\vec{t}}_n \, \mathrm{d}\Gamma =
	0 & \forall \tensorTwo{v} \in \boldsymbol{\mathcal{Q}}, 
	\label{eq:weak_form_mom} \\
	R_s = &
	-\int_\Omega \psi  \Delta m_{w}   \, \mathrm{d}\Omega 
	 + {\Delta t} \sum_{f \notin \Gamma_f^{N}} \llbracket \psi \rrbracket F_{w,n}^f 
	 + {\Delta t} \int_\Omega \psi q_{w,n} \, \mathrm{d}\Omega = 
	0 &  \forall \psi \in \mathcal{P}, 
	\label{eq:weak_form_massW} \\
	R_{p} =&
	-\int_\Omega \chi  \Delta m_{o}   \, \mathrm{d}\Omega 
	+ \Delta t \sum_{f \notin \Gamma_f^{N}} \llbracket \chi \rrbracket F_{o,n}^f 
	 + \Delta t \int_\Omega \chi q_{o,n} \, \mathrm{d}\Omega = 
	0 & \forall \chi \in \mathcal{P}, 
	\label{eq:weak_form_massNW} 
	\end{align}
\end{subequations}
where $\{ \tensorTwo{v}, \psi, \chi \}$ are discrete test functions. The symbol $\llbracket \cdot \rrbracket$ denotes the jump of a quantity across face $f$ in $\mathcal{T}^h$.   For an internal face, $\llbracket \chi \rrbracket = ( {\chi}_{\left| L \right.} - {\chi}_{\left| K \right.} )$, with ${\chi}_{\left| L \right.}$ and ${\chi}_{\left| K \right.}$ the restriction of $\chi$ on cells $K$ and $L$ sharing $f$, respectively.   For a face belonging to the domain boundary, the jump expression reduces to $\llbracket \chi \rrbracket = - {\chi}_{\left| K \right.}$. The term $F_{\ell,n}^f$ denotes a discrete mass flux, i.e. $F_{\ell,n}^f \approx \int_f \tensorTwo{w}_{\ell,n} \cdot \tensorTwo{n}_f \, dA$.  These are computed using a standard two-point flux approximation scheme, with upwinding of the phase density and mobility \cite{aziz1979}:
\begin{align}
  F_{\ell,n}^f &= - \rho^{\texttt{upw}}_{\ell,n} \lambda^{\texttt{upw}}_{\ell,n} \Upsilon^f \left( \llbracket p_n \rrbracket + \rho^f_{\ell,n} g \llbracket z \rrbracket \right). 
  \label{eq:num_flux_approx}
\end{align}
Here, $\Upsilon^f$ is the transmissibility coefficient for the face, which is computed knowing the mesh geometry and permeability.  The jump in the elevation datum should be understood as the difference in the $z-$coordinate of the respective cell centroids.

The unknown fields are interpolated as
\begin{subequations}
\begin{align}
\vec{u}_n(\vec{x}) &= \sum_{i=1}^{n_{\vec{u}}} \tensorTwo{\eta}_i(\vec{x}) u_{i,n}, \label{eq:U_approx}\\
s_n(\vec{x}) & = \sum_{j=1}^{n_s} \varphi_j(\vec{x}) s_{j,n}, \label{eq:S_approx}\\
p_n(\vec{x}) &= \sum_{k=1}^{n_p} \varphi_k(\vec{x}) p_{k,n}, \label{eq:P_approx}
\end{align}
\end{subequations}
\noindent
with $\{ \tensorTwo{\eta}_i \}$ and $\{ \varphi_j \}$ bases for $\boldsymbol{\mathcal{Q}}$ and $\mathcal{P}$, respectively.  $\{u_{i,n}\}$ are nodal values of the displacement components, while $\{s_{k,n}\}$ and $\{p_{k,n}\}$ are cell-centered values for the saturation and pressure fields.  An identical basis is introduced for the test functions. 

The fully discrete system of equations at time $t_n$ is then obtained by introducing these bases into the weak form Eqs.~\eqref{eq:weak_form_mom}-\eqref{eq:weak_form_massNW} and applying the standard  finite element procedure.  This leads to a set of algebraic equations for the unknown degrees-of-freedom $\{u_{i,n}\}$, $\{s_{j,n}\}$ and $\{p_{k,n}\}$. These degrees-of-freedom are assembled in an algebraic vector $\blkVec{x}_n = \{\Vec{u}_n$, $\Vec{s}_n$, $\Vec{p}_n\}$. The nonlinear system of equations is assembled in a residual vector 
\begin{align} 
\blkVec{r}_n \; (\blkVec{x}_n,\blkVec{x}_{n-1}) = 
\begin{bmatrix}
\Vec{r}_n^{u\;}  \\
\Vec{r}_n^{s\;}	 \\
\Vec{r}_n^{p\;} 
\end{bmatrix} = \blkVec{0}.
\label{eq:NL_res}
\end{align}
The nonlinearity of the system here results from nonlinear constitutive behavior embedded in the relative permeability relationship $k_{r\ell} (s)$. Furthermore, this formulation is general enough to accommodate other nonlinear constitutive relationships for the various solid and fluid components. A Newton iteration scheme is used to solve this system, which requires the linearization of the three-field problem. The linearized problem is defined by a Jacobian system with a 3 $\times$ 3 block structure of the form
\begin{align}
\begin{bmatrix}
\Mat{A}\sub{uu} & \Mat{A}\sub{us} & \Mat{A}\sub{up} \\
\Mat{A}\sub{su} & \Mat{A}\sub{ss} & \Mat{A}\sub{sp} \\
\Mat{A}\sub{pu} & \Mat{A}\sub{ps} & \Mat{A}\sub{pp}
\end{bmatrix}^k
\begin{bmatrix}
\delta \Vec{u}\\ \delta \Vec{s} \\ \delta \Vec{p}
\end{bmatrix}
&=
-
\begin{bmatrix}
\Vec{r}_n^{u} \\ \Vec{r}_n^s \\ \Vec{r}_n^{p}
\end{bmatrix}^k,
\label{eq:jac_system}
\end{align}
where $\blkMat{A}^k = \partial \blkVec{r}^k / \partial \blkVec{x}_n$ is the Jacobian matrix, with $\delta \Vec{u}$, $\delta \Vec{s}$, and $ \delta \Vec{p}$ the Newton search directions for each field. The superscript $k$ stands for the Newton iteration count. Full expressions for each elemental sub-vector of $\blkVec{r}$ and sub-matrix of $\blkMat{A}$ are reported in \cite{white2018}. We note that the solution of this linear system is typically the most expensive component of a fully-implicit code, and good solver performance is therefore essential.

\section{Incompressibility}
\label{sec:incomp}

It may not be immediately apparent that the system \eqref{eq:jac_system} may be subject to an inf-sup condition on its solvability.  Indeed, for many problem configurations the discrete system is perfectly well-posed.  Instabilities can arise, however, when two conditions are satisfied:
\begin{enumerate}
\item during undrained loading, i.e. as either $ \tensorTwo{\kappa} \to 0$ or $ \Delta t \to 0$;
\item when the solid and fluid phases approach incompressibility, i.e. $K_s \to \infty \text{ and } K_{\ell} \to \infty$ for $\ell=\{w,o\}$.  
\end{enumerate}
Note that it is sufficient to merely approach these limits, a situation that occurs frequently in practice.  This is particularly true at early simulation times, when small time steps $\Delta t$ are often required to resolve rapidly evolving solution fields.  Liquid and solid compressibilities are also often small for many geologic systems.

To highlight the origin of difficulties, we first revisit the continuum governing equations assuming the conditions above are exactly satisfied.   In this case, several relationships simplify, in particular
\begin{subequations}
\begin{align}
&b=1, \\
&\dot \phi = \nabla \cdot \dot{\tensorTwo{u}}, \\
&\dot \rho_s = \dot \rho_w = \dot \rho_o = 0 ,\\
&\tensorTwo{w}_w  = \tensorTwo{w}_o =  \tensorTwo{0} ,\\
&q_w  = q_o = 0 .
\end{align}
\end{subequations} 
We set $q_\ell=0$ here under the assumption that these source terms represent wells, which cannot inject when permeability goes to zero.  The mass balance equations \eqref{eq:massBalanceW_S}-\eqref{eq:massBalanceNW_S} reduce to
\begin{subequations}
	\begin{align}
	&\frac{\partial}{\partial t} \left({\phi s} \right)= 0 \\
	%
	&\frac{\partial}{\partial t} \left({\phi (1-s)} \right)= 0 
	\end{align}
\end{subequations}
Adding these two equations implies $\dot \phi = 0$, and therefore $\dot s = 0$.  The reduced system of governing equations is therefore
\begin{subequations}\label{eq:IBVP_reduced}
	\begin{align}
	&\nabla \cdot \tensorTwo{\sigma'} - \nabla p +  \rho\vec{g} = \vec{0} &\qquad \mbox{ in } \Omega \times \mathcal{I},\\
	%
	&\nabla \cdot \dot{\tensorTwo{u}} = 0 &\qquad \mbox{ in } \Omega \times \mathcal{I},
	\end{align}
\end{subequations}
with $s(\tensorTwo{x},t)=s_0(\tensorTwo{x})$.  We observe that under these conditions the solid deformation field must satisfy a divergence constraint condition, while the saturation field becomes fixed in time at its initial conditions.  This result is physically intuitive.  If the fluids can neither flow nor compress, they will not allow the solid skeleton to deform volumetrically, nor is there a mechanism for saturations to evolve.  The result is a two-field problem only in displacement and pressure.  With some additional manipulations one can show that this set of governing equations is equivalent to Stokes' equations.

It is also instructive to perform the same exercise for the algebraic system \eqref{eq:jac_system}.  When phase compressibility is zero and undrained conditions are reached, the Jacobian matrix becomes
\begin{equation}
\blkMat{A^*} = \begin{bmatrix}
\Mat{A}\sub{uu}^* &  & \Mat{A}\sub{up}^* \\
\Mat{A}\sub{su}^* & \Mat{A}\sub{ss}^* &  \\
\Mat{A}\sub{pu}^* & \Mat{A}\sub{ps}^* & 
\end{bmatrix}
\end{equation}
with the following expressions for the individual blocks:
\begin{subequations}
\begin{align}
\left[ \Mat{A}\sub{uu}^* \right]\sub{ij} &= \int_\Omega \nabla^s \tensorTwo{\eta}^i : \mathbb{C} : \nabla^s \tensorTwo{\eta}^j dV \\
\left[ \Mat{A}\sub{up}^* \right]\sub{ij} &= - \int_\Omega \nabla \cdot \tensorTwo{\eta}^i \varphi^j dV\\
\left[ \Mat{A}\sub{su}^* \right]\sub{ij} &= -\rho_w \int_\Omega s_n \varphi^i \nabla \cdot \tensorTwo{\eta}^jdV \\
\left[ \Mat{A}\sub{pu}^* \right]\sub{ij} &= -\rho_o \int_\Omega (1-s_n) \varphi^i \nabla \cdot \tensorTwo{\eta}^jdV \\
\left[ \Mat{A}\sub{ss}^* \right]\sub{ij} &= -\rho_w \int_\Omega \phi  \varphi^i \varphi^j dV \\
\left[ \Mat{A}\sub{ps}^* \right]\sub{ij} &= \rho_o \int_\Omega \phi  \varphi^i \varphi^j dV
\end{align}
\end{subequations}
Now imagine that we add the third block row to the second, scaled by their (constant) densities.  We observe that
\begin{subequations}
\begin{align} 
\frac{1}{\rho_w} \Mat{A}\sub{su}^* + \frac{1}{\rho_o} \Mat{A}\sub{pu}^* &= \Mat{A}\sub{up}^{T*} \\
\frac{1}{\rho_w} \Mat{A}\sub{ss}^* + \frac{1}{\rho_o} \Mat{A}\sub{ps}^* &= \Mat{0} 
\end{align}
\end{subequations}
We may therefore identify a sub-system for displacements and pressures that is uncoupled from the saturation field,
\begin{equation}
\blkMat{B^*} = \begin{bmatrix}
\Mat{A}\sub{uu}^* &  \Mat{A}\sub{up}^* \\
\Mat{A}\sub{up}^{T*} &  
\end{bmatrix}
\end{equation}
One would arrive at the same system through a direct discretization of the reduced governing equations above.  It is clear that this matrix is in saddle-point form, and that the spaces chosen for the pressure and displacement fields must satisfy an inf-sup compatibility condition to ensure $\blkMat{B^*}$ has full rank.  For our chosen interpolation, however, this is not the case.  As the incompressible limit is approached, $\blkMat{B^*}$---and equivalently, $\blkMat{A^*}$---may contain near-singular modes that can express as spurious oscillations in the pressure field.

\section{Stabilized Formulation}
\label{sec:stabilized}

As a fix for this difficulty, we propose a simple modification to the way the discrete fluxes are treated in Eqs. \eqref{eq:weak_form_massW}--\eqref{eq:weak_form_massNW}.  As described earlier, the mesh is decomposed into macroelements.  Let $\Gamma_M$ denotes the union of all faces $f$ that lie interior to any macroelement.  That is, any two cells connected across a face $f \in \Gamma_M$ are members of the same parent macroelement.

For any face in $\Gamma_M$, we augment the physical flux with an additional stabilization flux $G^f_\ell$.  For a given time increment $\Delta t$, we replace
\begin{equation}
\Delta t F_{\ell,n}^f  \gets \Delta t F_{\ell,n}^f + G^f_\ell \qquad \forall f \in \Gamma_M
\end{equation}
where for each phase the stabilization flux is a function of an inter-element jump $ \,\llbracket \Delta p \rrbracket$ across the face, scaled by particular constants,
\begin{subequations}\label{eq:artflux}
\begin{align}
G_w^f &= - \alpha^f_w  \,\llbracket \Delta p \rrbracket \quad \text{with} \quad  \alpha_w^f = \tau  \, V^e  \left[ \rho_{w} s\right]^\texttt{upw}_{n-1},\\
G_o^f & = - \alpha^f_o \, \llbracket \Delta p \rrbracket \quad \text{with} \quad  \alpha_o^f = \tau   \,V^e  \left[ \rho_{o} (1-s) \right]^\texttt{upw}_{n-1}.
\end{align}
\end{subequations}
These artificial fluxes will be used to control spurious pressure modes associated with non-physical pressure jumps across faces.  Here, $\tau$ is a stabilization parameter and $V_e$ is the average volume of the child elements in the macroelement.  The remaining terms are the upwinded density and phase saturation for the respective phases.  Note that these are lagged in time to simplify the linearization, as a lagged approximation of these quantities is usually sufficient for stabilization purposes.  We will discuss the choice of stabilization constant below, which is critical to success.

This flux form is quite similar to the physical flux computation Eq.~\eqref{eq:num_flux_approx}, so  it may be readily added to an existing face-based assembly loop.  Any addition of artificial fluxes, however, will break the element-wise mass conservation property of the underlying finite volume scheme.  Because these artificial fluxes are only added to internal faces of the macroelement, however, exact mass conservation is still preserved on the macroelement level.
 
 When assembled, these flux terms add additional entries to two blocks of the system matrix,
 \begin{equation}
 \blkMat{A} = 
 \begin{bmatrix}
\Mat{A}\sub{uu} & \Mat{A}\sub{us} & \Mat{A}\sub{up} \\
\Mat{A}\sub{su} & \Mat{A}\sub{ss} & \Mat{A}\sub{sp}+\Mat{C}\sub{sp} \\
\Mat{A}\sub{pu} & \Mat{A}\sub{ps} & \Mat{A}\sub{pp}+\Mat{C}\sub{pp}
\end{bmatrix},
\end{equation}
where the stabilizing entries are assembled face-wise for any $f \in \Gamma_M$ as
\begin{subequations}
\begin{align}
\left[\Mat{C}\sub{sp} \right]^f\sub{ij} & = -\alpha_w^f \, \llbracket \varphi^i \rrbracket \llbracket \varphi^j \rrbracket ,\\
\left[\Mat{C}\sub{pp} \right]^f\sub{ij} & = -\alpha_o^f \, \llbracket \varphi^i \rrbracket \llbracket \varphi^j \rrbracket .
\end{align}
\end{subequations}
In the incompressible limit, these contributions will not vanish, so that
\begin{equation}
\blkMat{A^*} = \begin{bmatrix}
\Mat{A}\sub{uu}^* &  & \Mat{A}\sub{up}^* \\
\Mat{A}\sub{su}^* & \Mat{A}\sub{ss}^* &  \Mat{C}\sub{sp}\\
\Mat{A}\sub{pu}^* & \Mat{A}\sub{ps}^* & \Mat{C}\sub{pp}
\end{bmatrix}.
\end{equation}

In practice, we always solve the three-field problem.  However, it is instructive to apply the same reduction procedure as before for the incompressible limit.  This leads to a reduced system
\begin{equation}
\blkMat{B^*} = \begin{bmatrix}
\Mat{A}\sub{uu}^* &  \Mat{A}\sub{up}^* \\
\Mat{A}\sub{up}^{T*} &  \Mat{C}
\end{bmatrix}
\end{equation}
where
\begin{equation}
\Mat{C} = \frac{1}{\rho_w} \Mat{C}\sub{sp} + \frac{1}{\rho_o} \Mat{C}\sub{pp}.
\end{equation}
Thus, the original saddle-point system is modified so that new entries appear in the lower-right-hand block.   For each face $f \in \Gamma_M$, this matrix contains contributions
\begin{equation}
\left[ \Mat{C} \right]\sub{ij}^f = - \tau V^e \, \llbracket \varphi^i \rrbracket \llbracket \varphi^j \rrbracket. 
\end{equation}
Because of the macroelement construction, the resulting matrix is extremely sparse.  In 3D, it is block-diagonal with one $8 \times 8$ block $\Mat{C}_M$ for each macroelement in the mesh, with entries
\begin{align}
\Mat{C}_M = \, - \tau V^e
			\begin{bmatrix}
			\ 3 & -1 & \  & -1 & -1 & \  & \  & \  \\
			-1 &  \ 3 & -1 & \  & \  & -1 & \  & \  \\
			\  & -1 & \ 3 & -1 & \  & \  & -1 & \  \\
			-1 & \  & -1 & \ 3 & \  & \  & \  & -1 \\
			-1 & \  & \  & \  & \ 3 & -1 &  & -1 \\
			\  & -1 & \  & \  & -1 & \ 3 & -1 & \  \\
			\  &  & -1 & \  & \  & -1 & \ 3 & -1 \\
			\  &  & \  & -1 & -1 & \  & -1 & \ 3
			\end{bmatrix},
\end{align}
where $V^e$ is the average volume of the child elements in the macroelement. Note that in our mesh geometry, $V^e=A^fd^f$ for any face interior to the macroelement, with $A^f$ the face area and $d^f$ the distance between the centroids of the neighboring cells. By construction, this pattern preserves the underlying two-point flux approximation (TPFA) stencil adopted in the original finite volume scheme, and will not cause any fill-in.  Note that the matrices $\Mat{C}\sub{sp}$ and $\Mat{C}\sub{pp}$ will have the same sparsity pattern as $\Mat{C}$, though their entries are weighted by the local saturation and densities at the faces.

\begin{figure}
\centering
\includegraphics[width=0.33\textwidth]{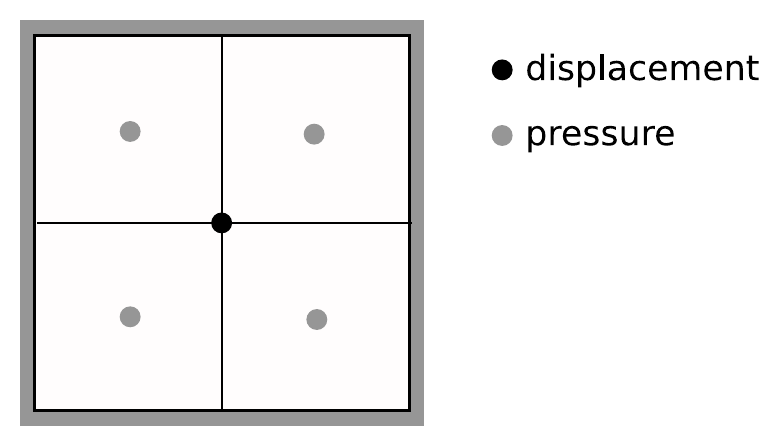}
\caption{Single macroelement patch test geometry in 2D, with one pressure in each cell and two displacement components at the central node.  The 3D macroelement is similar, involving eight pressures and three displacement components.}
\label{fig:patchtest}
\end{figure}

For the reduced system, the stabilization contribution $\Mat{C}$ has a very similar form to the Local Pressure Jump (LPJ) stabilization as originally formulated for solving the Stokes equation \cite{silvester1990}. Indeed, this basic idea of using inter-element pressure jumps to control spurious modes inspired the method proposed here.  There are two key differences for the multiphase poromechanics application, however: 
\begin{enumerate}
\item In the three-field formulation, a separate contribution is made to each mass balance equation, weighted by appropriate phase density and saturation;
\item The optimal stabilization constant $\tau$ will differ due to the nature of the underlying equations.
\end{enumerate}
We still have to address this question of what is an appropriate value for the stabilization parameter $\tau$. As proposed in \cite{silvester1994}, good candidates for $\tau$ can be determined by examining the spectrum of the Schur complement matrix,
\begin{equation}
\Mat{S} =   \Mat{A}\sub{up}^{*T}  \Mat{A}\sub{uu}^{*-1} \Mat{A}\sub{up}^* - \Mat{C},
\end{equation}
which corresponds to a further block reduction of $\blkMat{B^*}$ to a pressure-only system.    We focus on a patch test involving a single macroelement, with rigid and impermeable boundary conditions (Figure~\ref{fig:patchtest}). If stability can be demonstrated for a single macroelement, theoretical results in \cite{boland1983,stenberg1984} prove stability for discretizations on arbitrary grids constructed by ``gluing together'' stable macroelements. 

The resulting Schur-complement is rank deficient without stabilization.  In 3D, there are eight pressure unknowns in the cells but only three displacement components at the central node.  Similar to \cite{elman2005}, the eigenvalues and eigenvectors of this system may be readily computed.  Let the cells have edge lengths $h_x$, $h_y$, and $h_z$. The element volume is $V=h_xh_yh_z$, and each face has area $A_{xy}=h_xh_y$, $A_{yz}=h_yh_z$, or $A_{xz}=h_xh_z$.  The resulting eigenvalues are
 \begin{subequations}
 \begin{align}
 &e_1 = 0\\
 &e_2 = e_3 = e_4 = 4  V \tau\\
 &e_5 = 6V  \tau  \\
 &e_6 = V \left( 2 \tau + \frac{9A_{xy}^2}{16\left(A_{xy}^2(\lambda+2G) + A_{xz}^2G+A_{yz}^2G\right)} \right)\\
 &e_7 = V \left( 2 \tau + \frac{9A_{xz}^2}{16\left(A_{xy}^2G + A_{xz}^2(\lambda+2G)+A_{yz}^2G\right)} \right)\\
  &e_8 = V \left( 2 \tau + \frac{9A_{yz}^2}{16\left(A_{xy}^2G + A_{xz}^2G+A_{yz}^2(\lambda+2G)\right)} \right)
  \end{align}
  \end{subequations}
 where $\lambda$ and $G$ are the two Lam\'e parameters characterizing the elastic mechanical response.
 
 When no stabilization is applied to the macroelement, that is when $\tau = 0$, five out of eight eigenvalues are zero. The null eigenvectors include one constant pressure mode (associated to $e_1$) and four spurious pressure modes (associated to $e_2$--$e_5$). The constant mode is expected here since the boundary conditions only determine the pressure solution up to an arbitrary constant.  We see that for $\tau > 0$, the stabilization will remove the spurious pressure modes from the null space of $\Mat{S}$.  Unfortunately, the choice of $\tau$ will also impact the physical modes associated with $e_6$--$e_8$.
 
  \begin{figure}[t]
 	\centering
 	\includegraphics[width=0.45\textwidth]{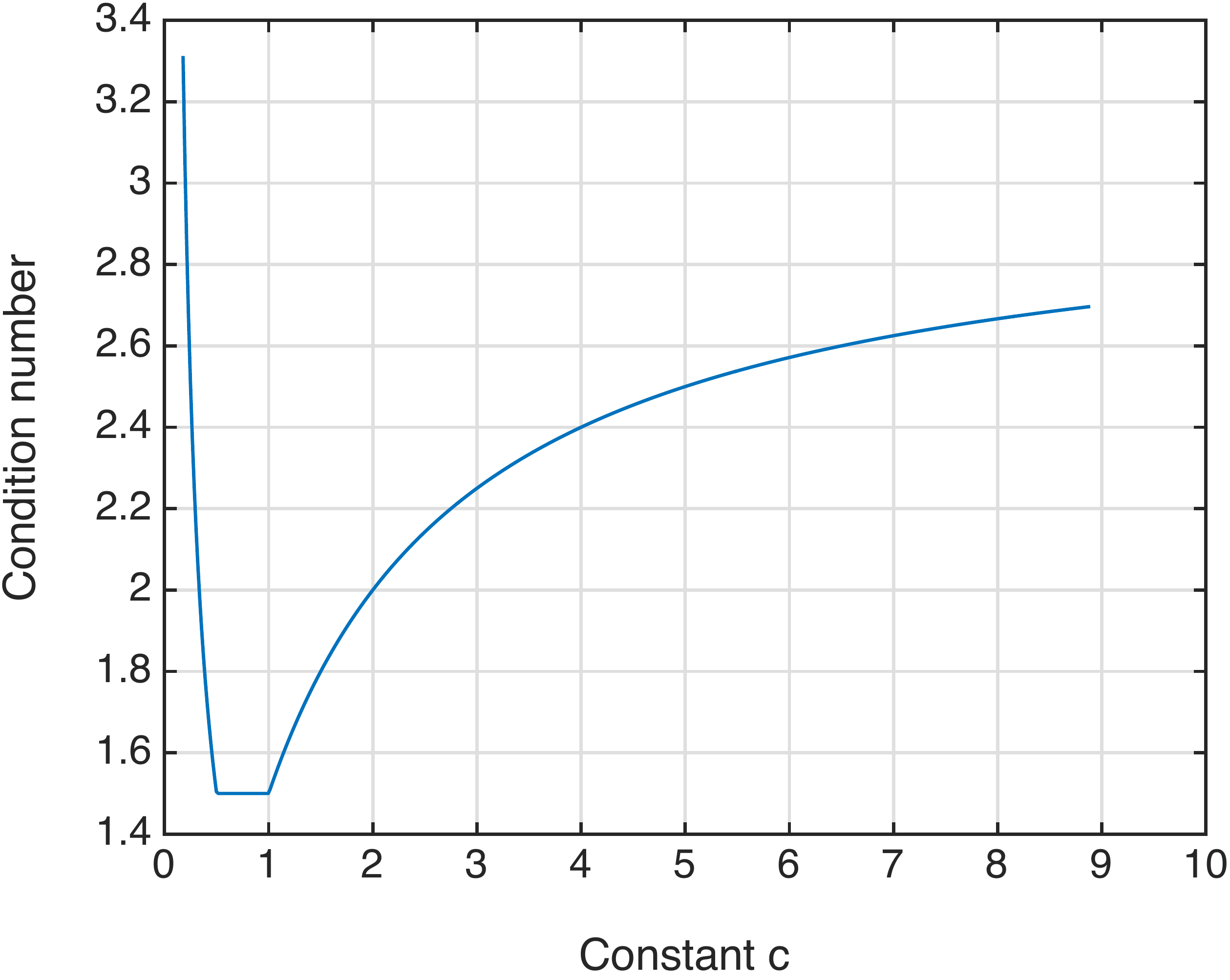}
 	\caption{Analytically computed condition number of  $\Mat{S}$  as a function of the stabilization parameter $\tau$ for the single macroelement patch test. Here, $c=\tau / \tau^*$ with $c=1$ being the recommended value.} \label{fig:min_cond_number}
 \end{figure}
 
\begin{table*}[t]
	\caption {Simulation parameters used in the three single-phase Barry-Mercer examples.} \label{tab:parameters} \small
	\begin{center}
		\begin{tabular}{ l  l  lll l}
			\toprule 
			Symbol & Parameter  & Drained & Undrained & Modified & Units \\
			\midrule			
			$\Delta t$  & Time step &   $6.14 \times 10^{-5}$ & $10^{-4}$ & $10^{-2}$& s  \\
			$T$  & Final time &   $1.54 \times 10^{-3}$ & $10^{-4}$ & $10^{-2}$ & s  \\
			$b$ & Biot coefficient &  1 & 1 & 1 & - \\
			$E$ & Young's modulus &  $10^{5}$ & $10^{5}$ & 2.5 & Pa \\
			$\nu$ & Poisson ratio &  0.1 & 0.1 & 0.25 & \\
			$k$ & Permeability &  $10^{-5}$ & $10^{-9}$ & $10^{-11}$ & $\text{m}^2$ \\
			$\mu_w$ & Viscosity &  $10^{-3}$  & $10^{-3} $ & $10^{-3}$ & Pa$\cdot$s \\
			$\rho_w$ & Density &  $10^{3}$ & $10^{3}$ & $10^{3}$ & $\text{kg/m}^3$\\
			\bottomrule 
		\end{tabular}
	\end{center}
\end{table*}


For stability, all eigenvalues must remain bounded away from zero except for the constant mode $e_1$.  Taking $\tau$ too large, however, will corrupt the physical solution and compromise the approximation. A good choice to balance these competing priorities is to choose a $\tau$ that minimizes the condition number $\kappa (\Mat{S}) = e_\text{max} / e_\text{min}$.  Considering a regular cube with equal sides $h_x = h_y = h_z = h$, the minimal condition number is retrieved for any stabilization parameter $\tau$ lying within the range 
  \begin{equation}
  \dfrac{9}{64 \left( \lambda + 4G \right)} \leq \tau \leq \dfrac{9}{32 \left( \lambda + 4G \right)}.  \label{eq:ideal_tau_range}
  \end{equation}
  Figure \ref{fig:min_cond_number} illustrates how the condition number $\kappa$ varies depending on the stabilization parameter. We can see that the condition number attains its minimal value of $\kappa = 3/2$ within the range prescribed in Eq.\eqref{eq:ideal_tau_range}.   Within this range, neither extremal eigenvalue actually depends on $\tau$.  Therefore, a reasonable choice for the stabilization parameter is 
  \begin{equation}
  \tau^* = \dfrac{9}{32 \left( \lambda + 4G \right)} . \label{eq:ideal_tau}
  \end{equation}
In order to explore the sensitivity of numerical solutions to this constant, it is convenient to present results in terms of the ratio $c=\tau / \tau^*$, i.e. the ratio of a given stabilization $\tau$ to the recommended value $\tau^*$, with $c=1$ being ``optimal''.

We emphasize that the recommended value here is based on the analysis of a fully-incompressible, homogeneous, cubic macroelement.  It therefore only depends on two elastic constants.  A more precise analysis may lead to a stabilization constant that additionally depends on the Biot coefficient, fluid compressibilities, macroelement heterogeneity, and so forth.  For example, when $b<1$ a better estimate would be $\tau = b^2 \tau^*$. 

In summary, the proposed stabilization method consists of adding the artificial flux terms in Eq.~\eqref{eq:artflux} to all macroelement-interior faces, weighted by the stabilization constant recommended in Eq.~\eqref{eq:ideal_tau}.  We now explore the effectiveness of this simple method with several numerical examples.

\section{Numerical Examples}
\label{sec:numerical_results}

We begin with a few single-phase examples ($s=1$) to demonstrate the performance of the method in a simpler setting.  We then conclude with a multiphase demonstration for a benchmark reservoir simulation problem.  These numerical experiments were implemented using \textit{Geocentric}, a simulation framework for computational geomechanics that relies heavily on finite element infrastructure from the the \textit{deal.ii} library \cite{bangerth2007}.

%
%
%

\subsection{Single-phase Examples}


\begin{figure}[p]
	\centering
	\includegraphics[width=0.4\textwidth]{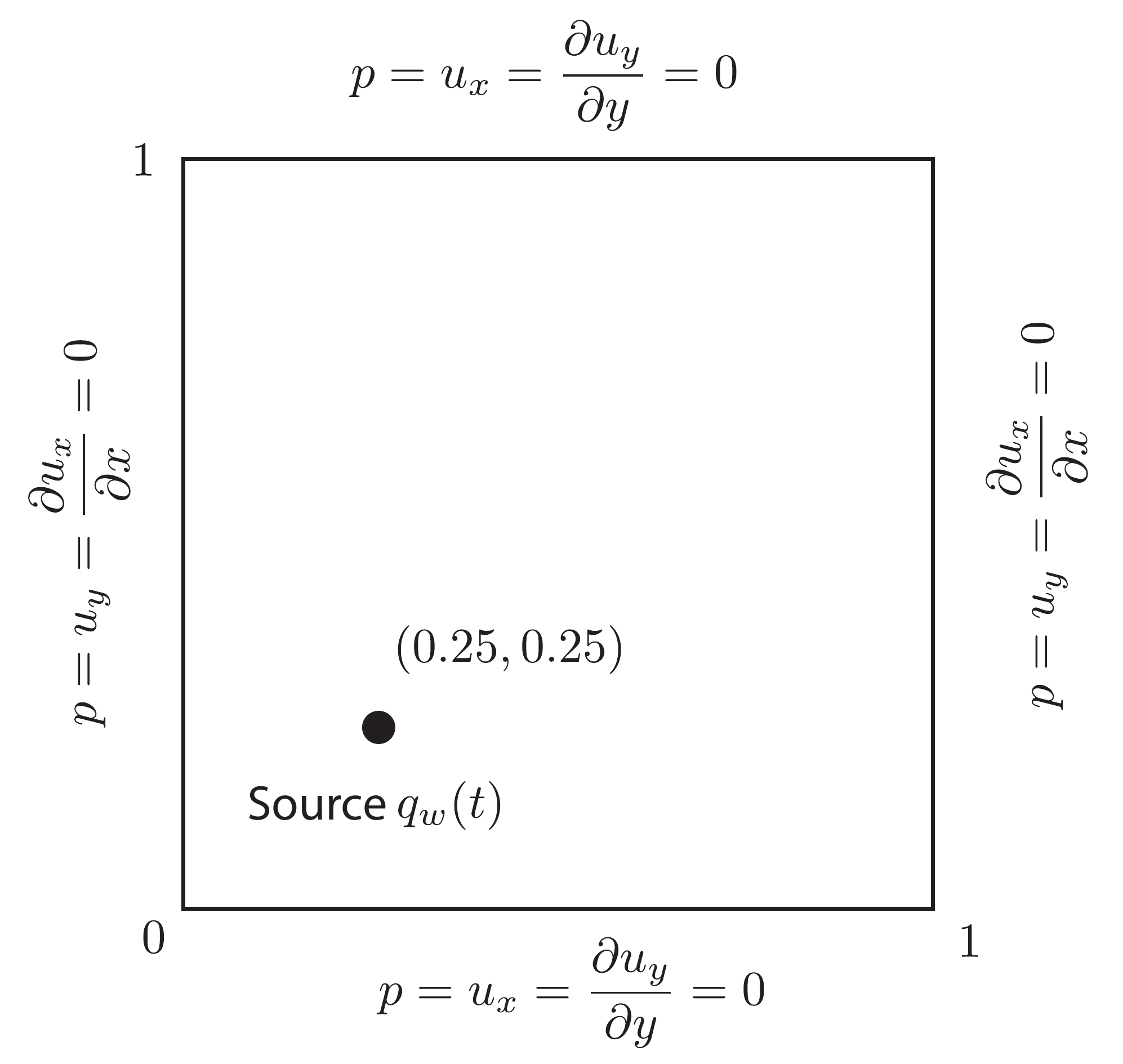}
	\caption{Geometry and boundary conditions for the Barry-Mercer problem.} \label{fig:bm_setup}
\end{figure}

 \begin{figure}[p]
 	\centering
 	\includegraphics[angle=90, width=0.5\textwidth]{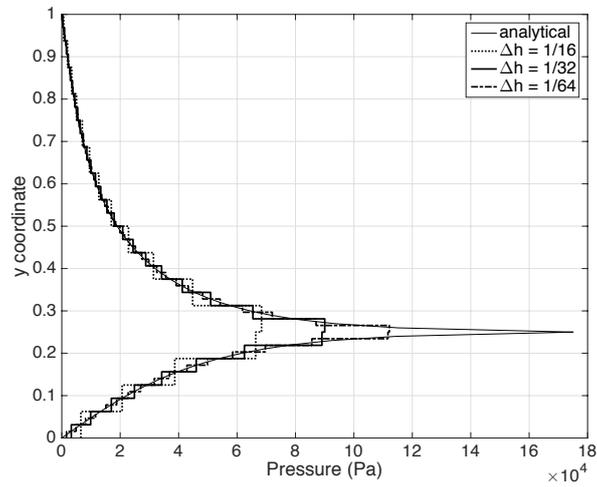}
 	\caption{Pressure plot along the vertical line $x=0.25$ for the drained Barry-Mercer problem at $\hat{t} = \pi / 2$ using the stabilized formulation.} \label{fig:bm_pressure2Dplot}
 \end{figure}
 
  \begin{figure}[t]
  	\centering
  	\includegraphics[width=0.5\textwidth]{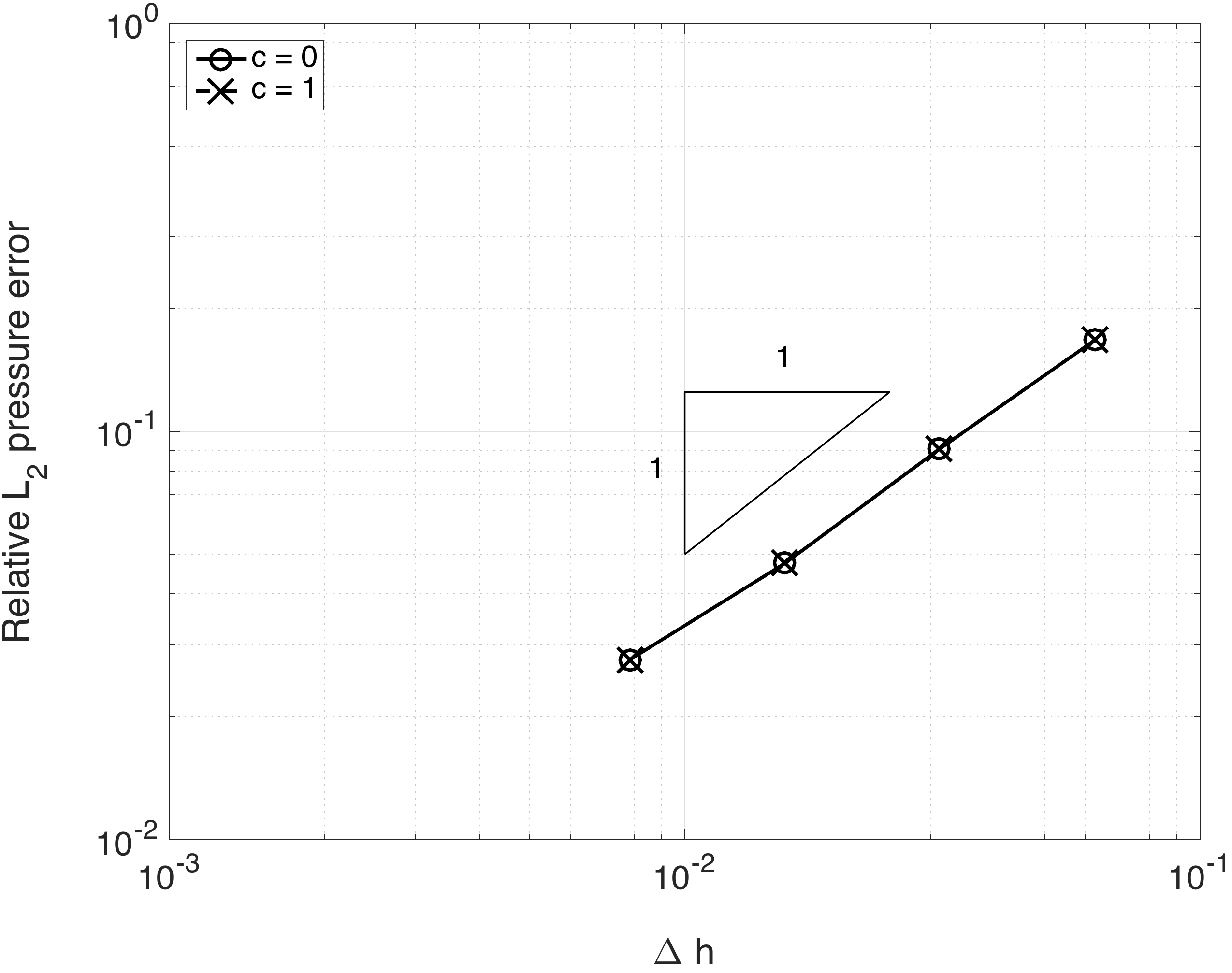}
  	\caption{Convergence of the relative $L_2$-error in pressure for the drained Barry-Mercer problem at $\hat{t} = \pi / 2$, using both the unstabilized ($c=0$) and stabilized ($c=1$) formulations.} \label{fig:bm_pressure_error}
  \end{figure}
  
  
  \begin{figure}[p]
  	\begin{center}
  		\begin{subfigure}[b]{0.45\textwidth}
  			\includegraphics[width=\textwidth]{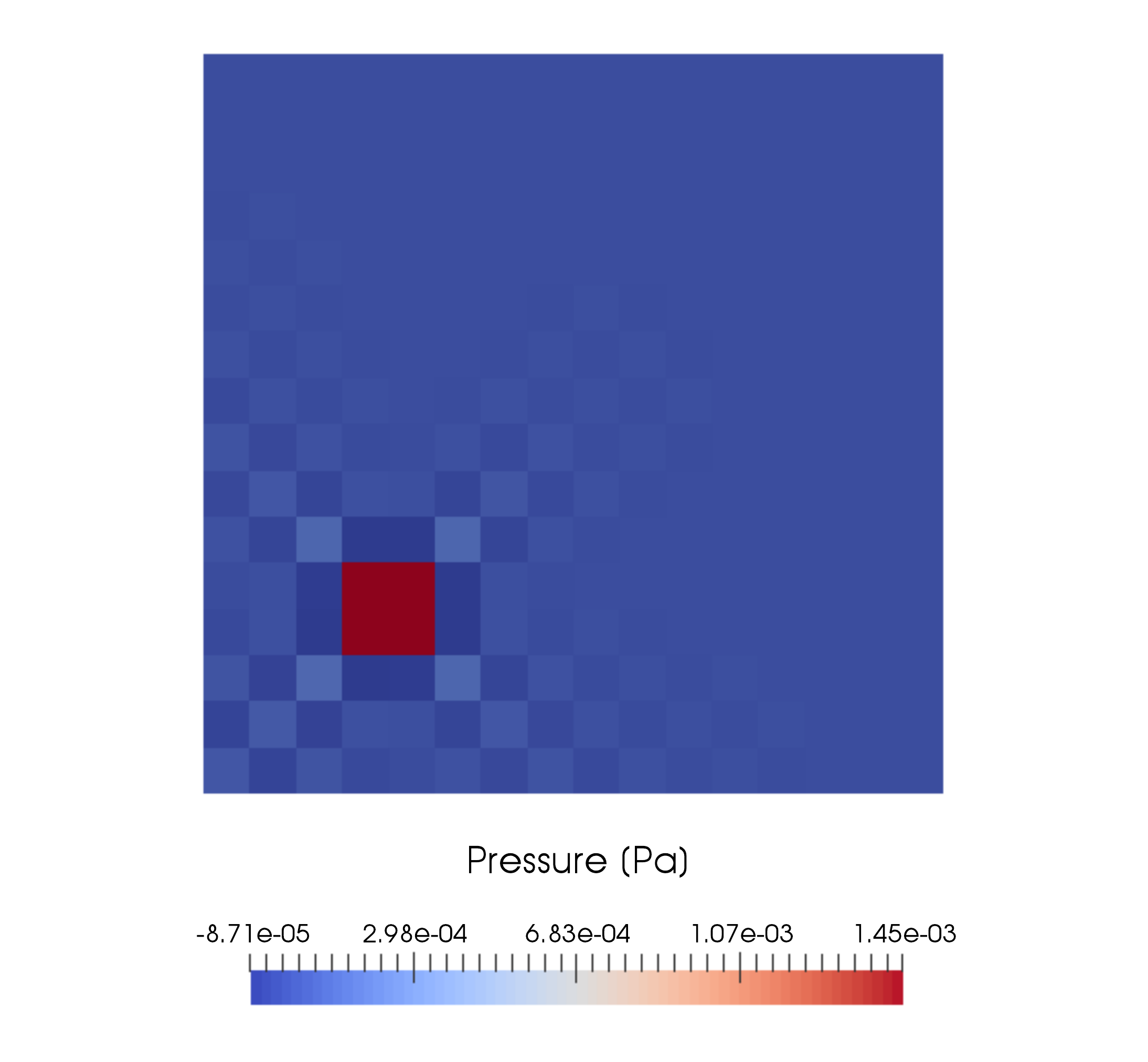}
  			\caption{Unstabilized scheme} \label{fig:bm_pressure_a}
  		\end{subfigure}
  		\begin{subfigure}[b]{0.45\textwidth}
  			\includegraphics[width=\textwidth]{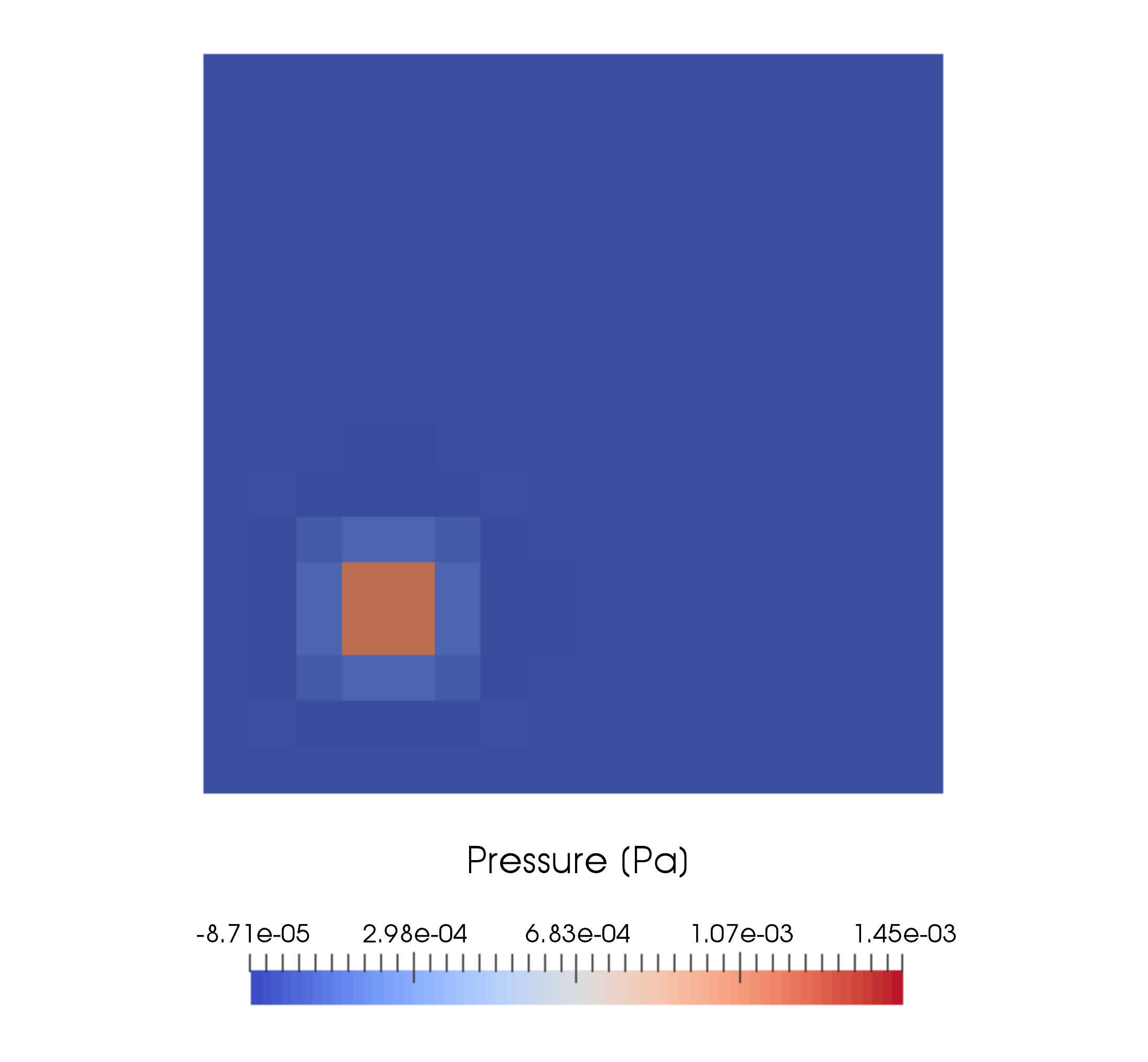}
  			\caption{Stabilized scheme (with $\tau = \tau^*$)} \label{fig:bm_pressure_b}
  		\end{subfigure}
  		\caption{Pressure distribution for the undrained Barry-Mercer problem.} \label{fig:bm_pressure}
  	\end{center}
  \end{figure}
  
    \begin{figure}[p]
  	\centering
  	\includegraphics[angle =90, width=0.45\textwidth]{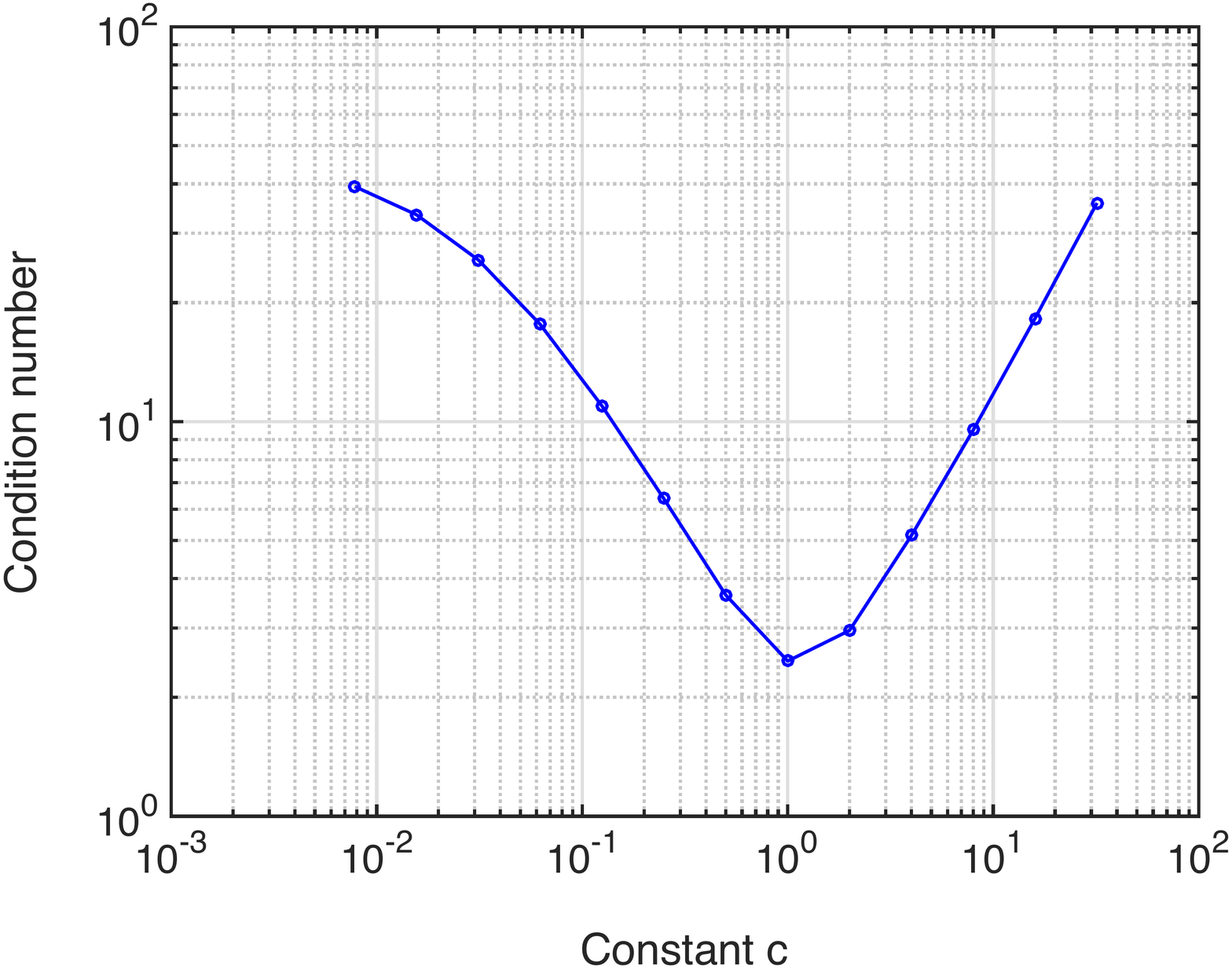}
  	\caption{Computed condition numbers of $\Mat{S'}$ for the undrained Barry-Mercer problem for various stabilization values $\tau = c \tau^*$.} \label{fig:bm2_cond_number}
  \end{figure}
  
    \begin{figure}[t]
  	\centering
  	\includegraphics[width=0.45\textwidth]{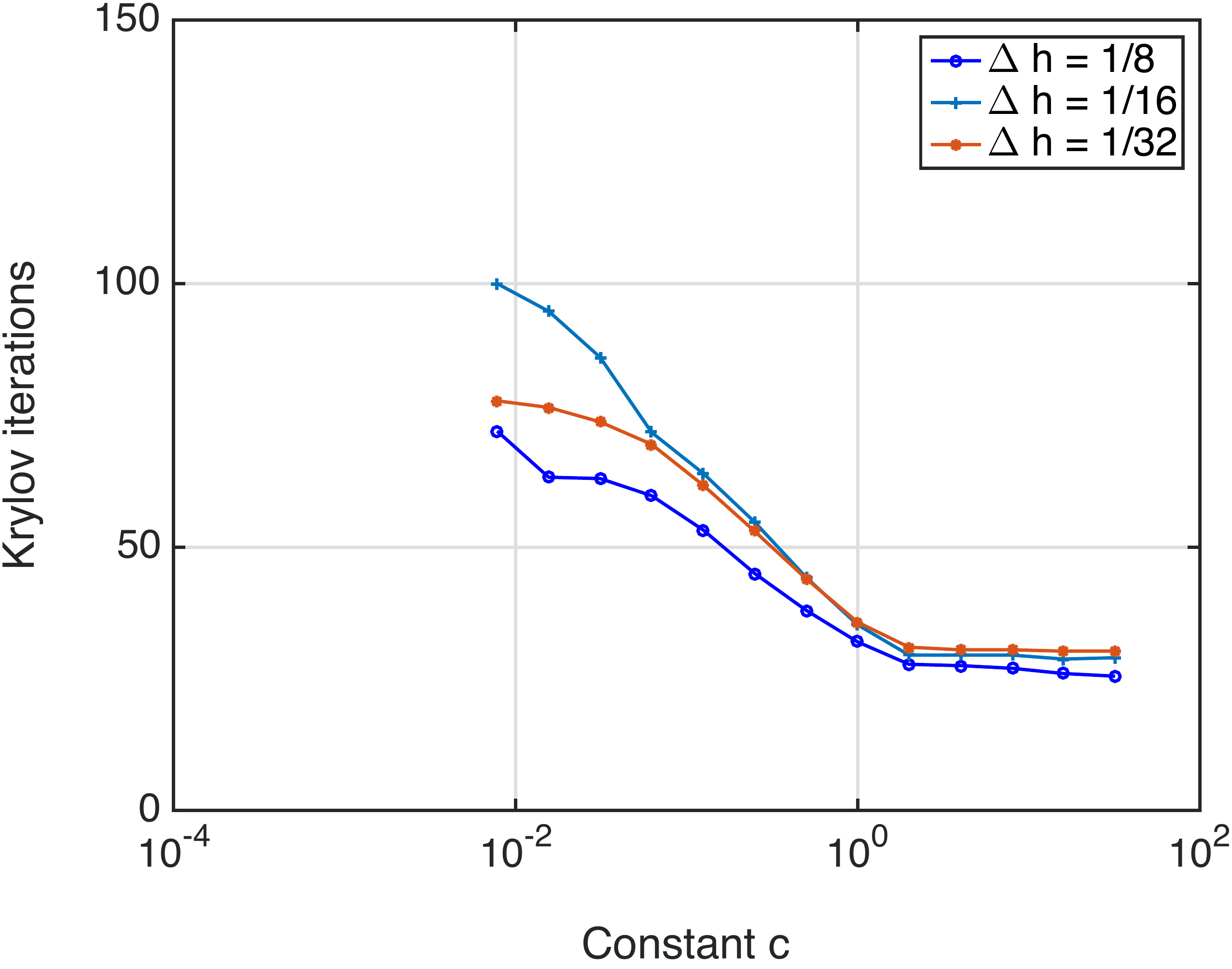}
  	\caption{Krylov iterations to convergence as a function of mesh refinement and stabilization constant $\tau = c \tau^*$ for the undrained Barry-Mercer problem.} \label{fig:bm02_krylov}
  \end{figure}

    
    \begin{figure}[p]
  	\begin{center}
  		\begin{subfigure}[b]{0.45\textwidth}
  			\includegraphics[width=\textwidth]{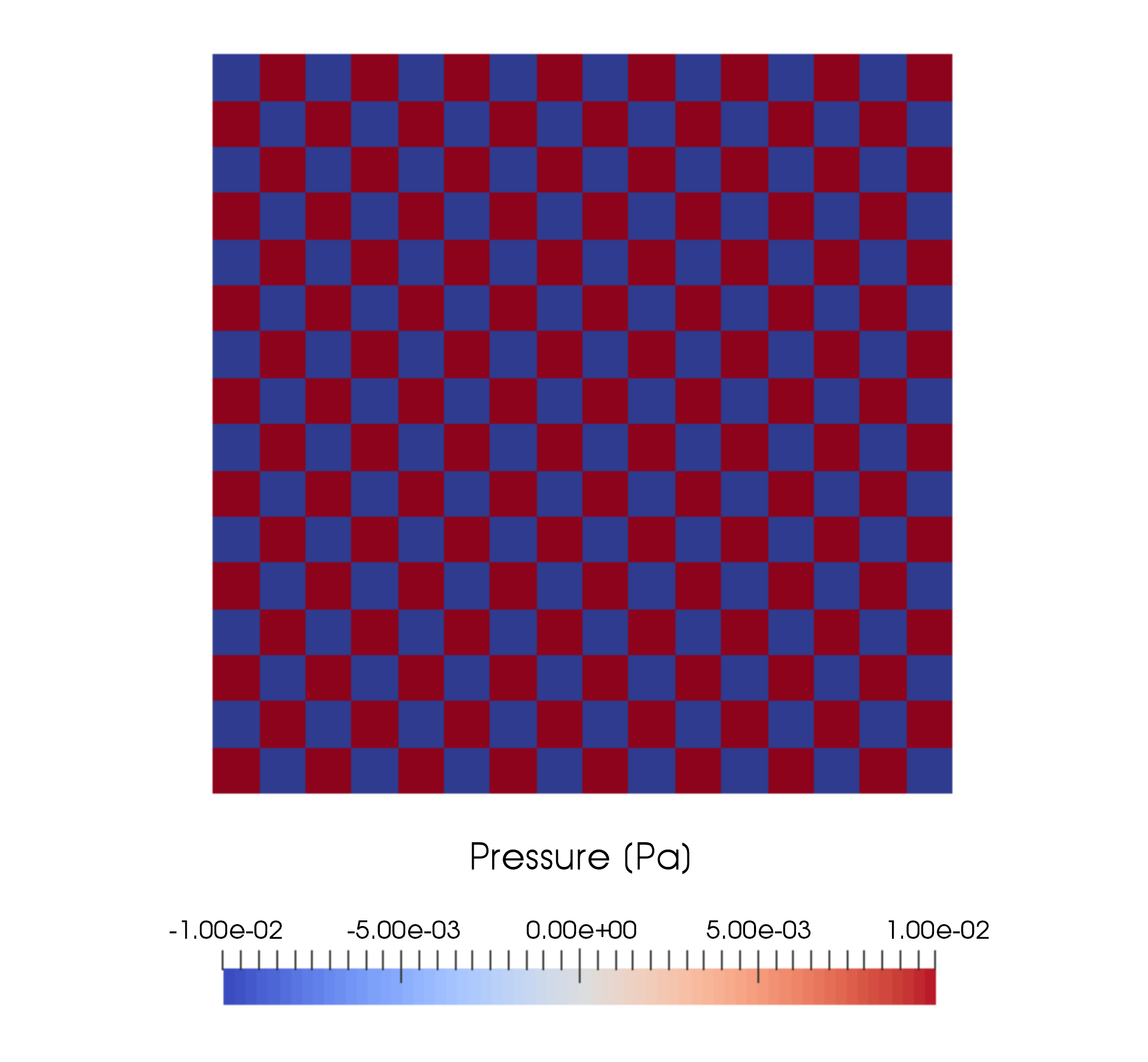}
  			\caption{Unstabilized scheme} \label{fig:bm3_pressure_a}
  		\end{subfigure}
  		\begin{subfigure}[b]{0.45\textwidth}
  			\includegraphics[width=\textwidth]{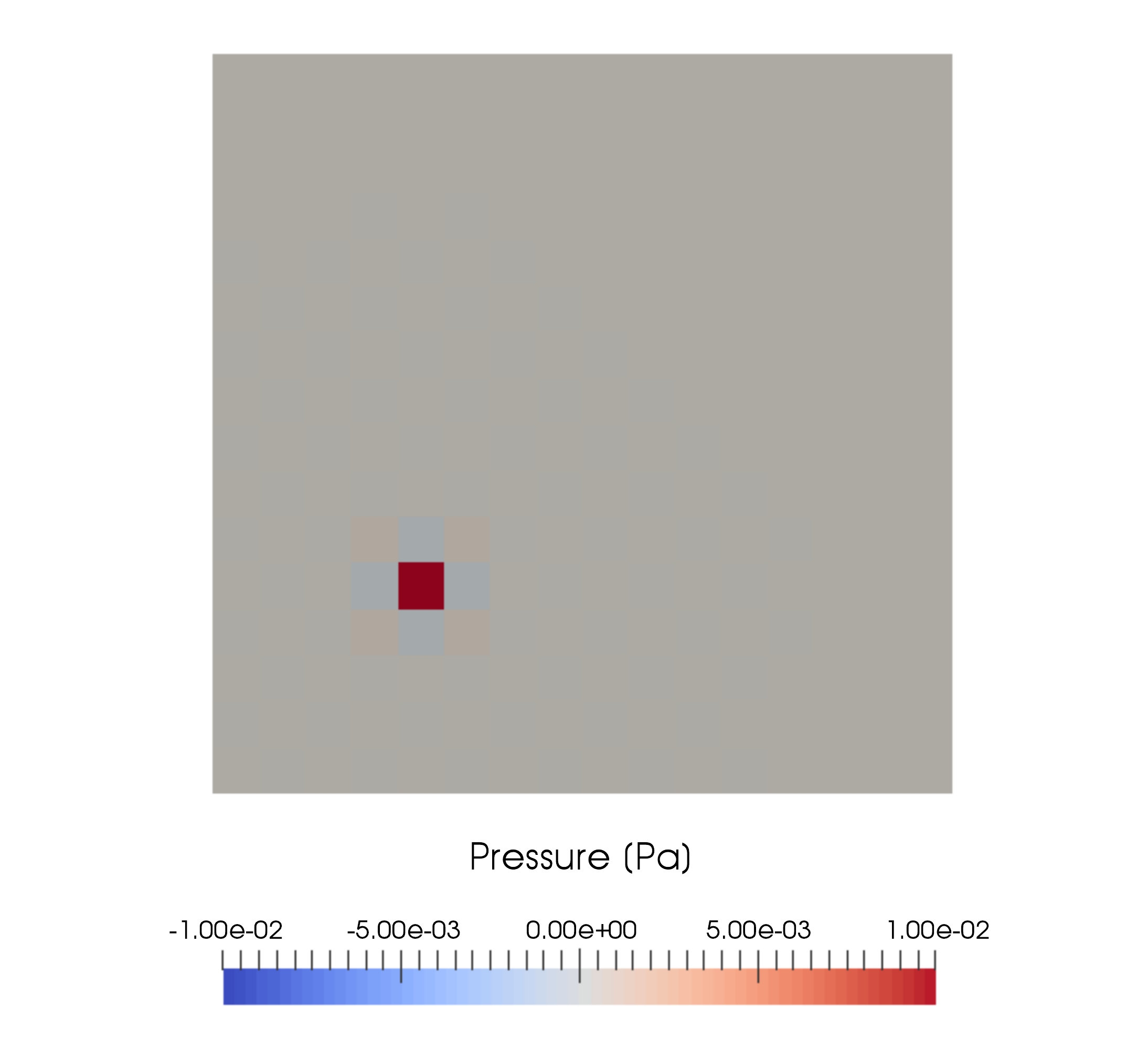}
  			\caption{Stabilized scheme (with $\tau = \tau^*$)} \label{fig:bm3_pressure_b}
  		\end{subfigure}
  		\caption{Pressure distribution for the modified undrained Barry and Mercer's problem.} \label{fig:bm3_pressure}
  	\end{center}
  \end{figure}
  
    \begin{table}[p]
		\caption{Extremal eigenvalues of the scaled Schur complement matrix for the modified undrained Barry-Mercer problem. Minimum and maximum eigenvalues as well as the condition number are presented as a function of mesh refinement. Two cases are considered: using (a) an unstabilized scheme and (b) the proposed stabilization.}
	\begin{subtable}[h]{1.0\textwidth}
		\caption{$\tau = 0$}
		\label{tab:bm03_evalue_0}
		\centering
			\begin{tabular}{ c  c c  c }
				\toprule
			\begin{tabular}{@{}c@{}}Number \\ of cells \end{tabular} & $e_{\text{min}}$ & $e_{\text{max}}$ & \begin{tabular}{@{}c@{}}Condition \\ number  ($e_{\text{max}}/e_{\text{min}}$) \end{tabular}\\
				\midrule		
			8 x 8  & 2.51 e$-03$  & 0.330 & 131.56  \\
			16 x 16  & 5.35 e$-04$  & 0.333 & 621.85  \\
			32 x 32  &  1.19 e$-04$ &  0.333  & 2799.72 \\
				\bottomrule
			\end{tabular}
	\end{subtable}
	
	\bigskip
	
	\begin{subtable}[h]{1.0\textwidth}
		\caption{$\tau = \tau^*$}
		\label{tab:bm03_evalue_1}
		\centering
		\begin{tabular}{ c  c c  c }
			\toprule
		\begin{tabular}{@{}c@{}}Number \\ of cells \end{tabular} & $e_{\text{min}}$ & $e_{\text{max}}$ & \begin{tabular}{@{}c@{}}Condition \\ number  ($e_{\text{max}}/e_{\text{min}}$) \end{tabular}\\
			\midrule		
		8 x 8  & 0.222  & 0.539 & 2.421   \\
		16 x 16  & 0.224  & 0.546  & 2.437  \\
		32 x 32  &  0.225  &  0.548  & 2.438 \\
			\bottomrule
		\end{tabular}
	\end{subtable}
	\label{tab:bm03_evalue}
\end{table}

    \begin{figure}[p]
  	\centering
  	\includegraphics[angle=90, width=0.45\textwidth]{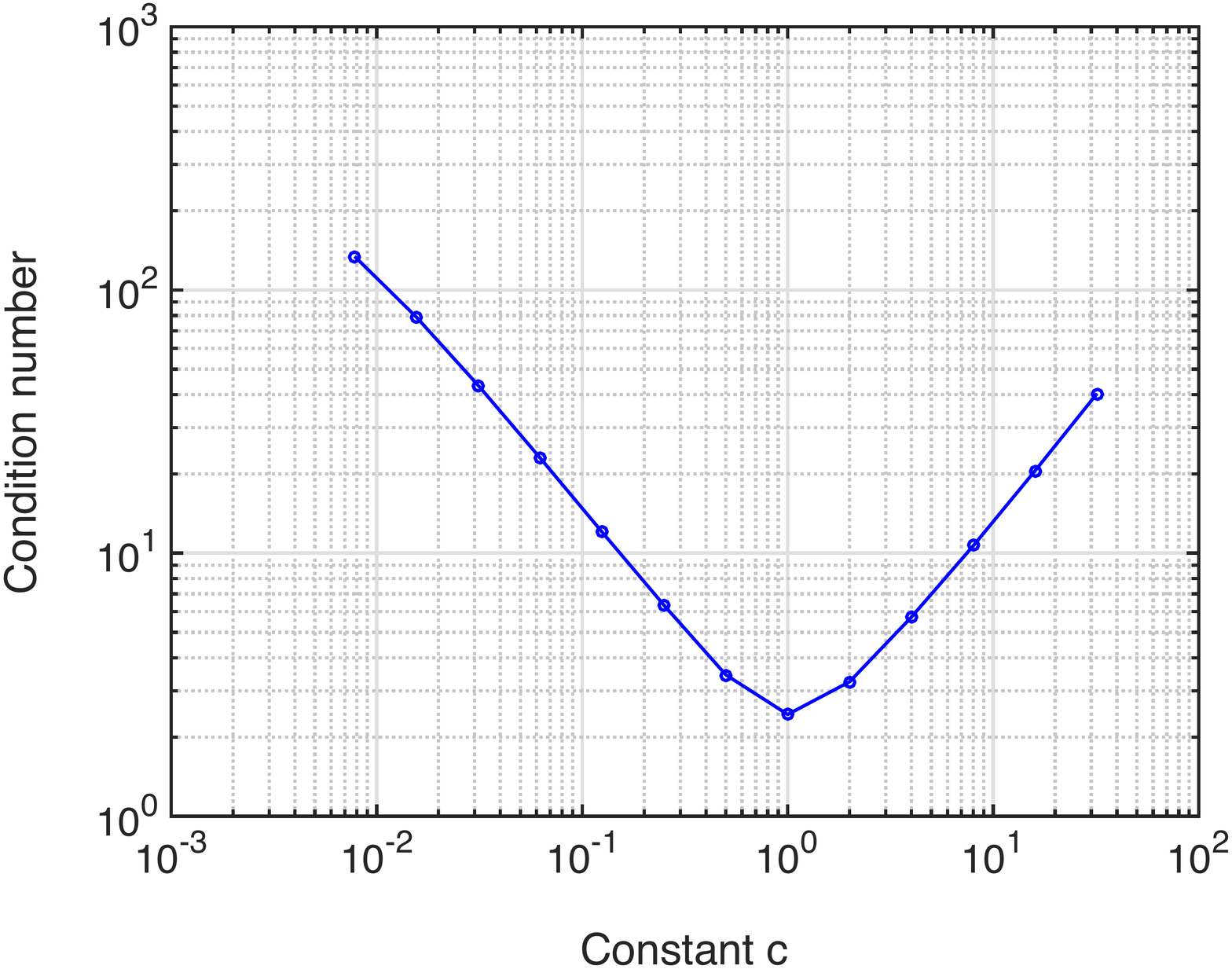}
  	\caption{Computed condition numbers of $\Mat{S'}$ for the modified undrained Barry-Mercer problem for various stabilization values $\tau = c \tau^*$.} \label{fig:bm03_cond_number}
  \end{figure}

  \begin{figure}[p]
  	\centering
  	\includegraphics[width=0.45\textwidth]{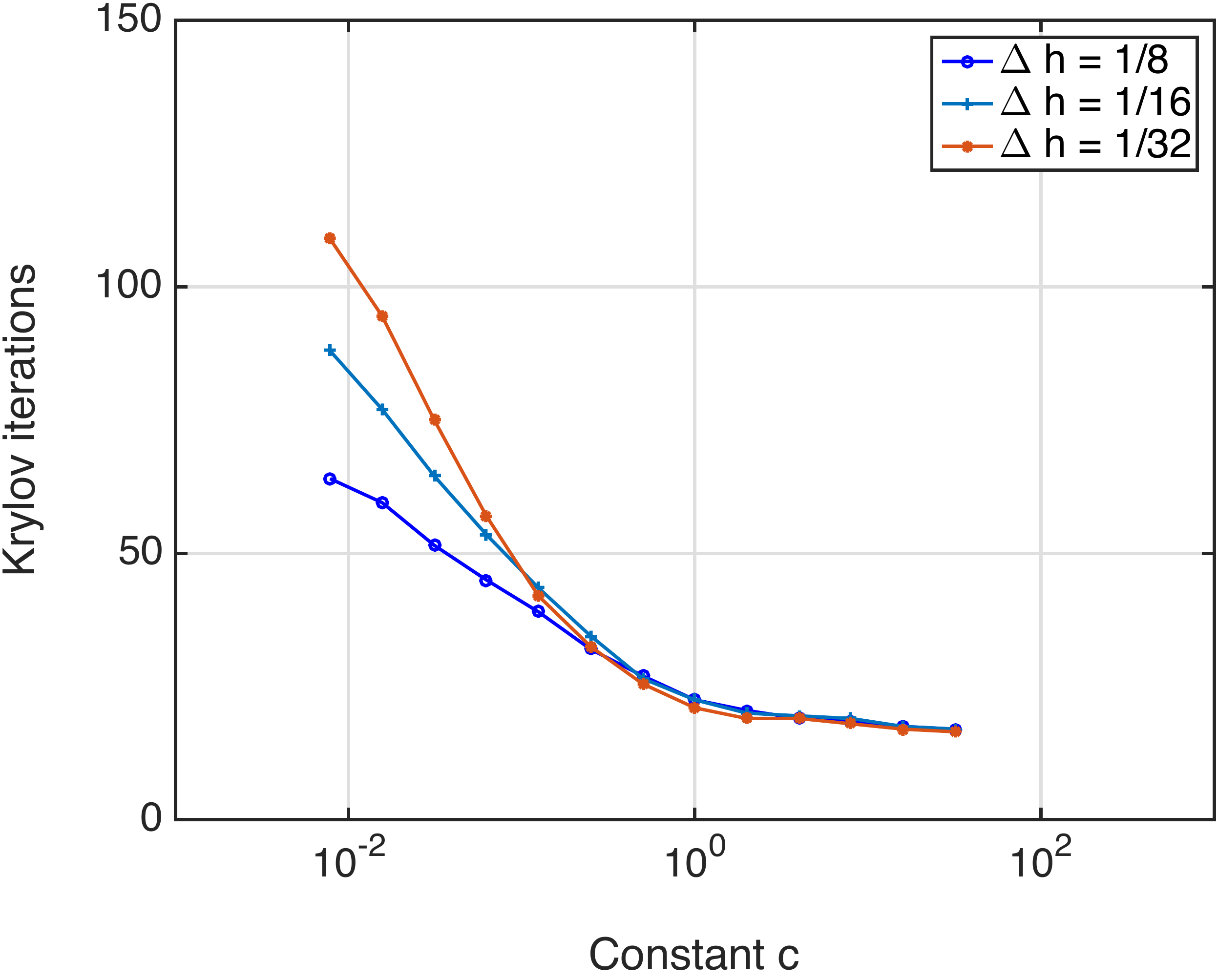}
  	\caption{Krylov iterations as a function of the mesh refinement and the stabilization constant for the modified undrained Barry Mercer problem.} \label{fig:bm03_krylov}
  \end{figure}
  
For the single-phase examples, we consider several variants of Barry and Mercer's problem for a two-dimensional, poroelastic medium \cite{barry1999}.  Parameter values are summarized in Table~\ref{tab:parameters}.   The first test illustrates that the proposed stabilization scheme does not compromise solution accuracy under drained conditions, when no instabilities are expected. Two subsequent examples confirm the effectiveness of the scheme in suppressing spurious pressure oscillations under undrained conditions.

 \subsubsection{Drained Barry-Mercer}
  \label{sec:bm_1}
 
Barry and Mercer \cite{barry1999} provide an analytical solution for a two dimensional problem. The problem setup consists of a square domain $\Omega = [0, 1] \times [0, 1]$ subjected to a periodic point source given as
\begin{equation}
q_{w}(t) = 2 \beta \, \delta \left(\vec{x} - \vec{x}_0 \right)\sin(\beta t) \qquad \text{with} \qquad \beta = \left( \lambda + 2 G \right) \frac{ \kappa }{  \mu }.
\end{equation}
Here, $\lambda$ and $G$ denote the two elastic Lam\'e parameters, while $\kappa$ and $\mu$ are the isotropic absolute permeability and viscosity, respectively. The point source is located at $\vec{x}_0 = \left( 0.25, 0.25 \right)$, and $\delta(\cdot)$ indicates a Dirac function. All sides of the computational domain are constrained with  zero pressure and zero tangential displacement boundary conditions, as depicted in Figure~\ref{fig:bm_setup}. The simulation parameters provided in Table \ref{tab:parameters} (drained conditions) are the same as those used in \cite{fu2019, boffi2016, rodrigo2016, phillips2007}. Note that the time step and final time correspond to a normalized time $\hat{t} = \beta \, t$ of $\Delta\hat{t} = 2 \pi / 100$ and $\hat{t} = \pi / 2$, respectively.
   
Figure~\ref{fig:bm_pressure2Dplot} shows the resulting pressure profile at the final time along a vertical line through the source point.  Both the analytical solution and the numerical solution for different mesh refinements are shown. Good agreement between the exact and computed results is also indicated by Figure~\ref{fig:bm_pressure_error}, which shows convergence behavior of the $L_2$-error for the pressure solution for both the stabilized and unstabilized formulations. One observes a linear and essentially identical error behavior for both, indicating that the macroelement stabilization does not compromise solution accuracy in regimes where it is not strictly needed. 
 
\subsubsection{Undrained Barry-Mercer}
\label{sec:bm_2}

The goal of this example is to show the effectiveness of proposed stabilization scheme in treating non-physical pressure oscillations. These spurious pressure oscillations appear in the limit of low permeability or fast loading rates. As in \cite{rodrigo2016,fu2019,boffi2016}, we use the same simulation parameters as the previous section, but we decrease the value of the permeability to $\kappa = 10^{-9} \, \text{m}^2$ and perform only one time step of $\Delta t = 10^{-4}$ s. 

Figure~\ref{fig:bm_pressure} shows the pressure contour plot for the uniformly discretized domain with $\Delta h = 1/16$. The pressure field exhibits mild oscillations close to the source-point. These oscillations are eliminated when using the proposed stabilization technique. 

It is also interesting to examine the conditioning of the stabilized system and its impact on iterative solver performance.  To do so, we define a scaled Schur-complement matrix,
\begin{equation}
\Mat{S'} = \Mat{Q}^{-1} \left(  \Mat{A}\sub{up}^T \Mat{A}^{-1}\sub{uu} \Mat{A}\sub{up} - \Mat{C} \right)
\end{equation}
where $\Mat{Q}$ is the mass matrix on the pressure space.  For the $\mathbb{P}_0$ space, this is a diagonal matrix with entries corresponding to the element volumes.  The inverse of this diagonal matrix introduces a volume scaling that allows eigenvalues to be properly compared across different mesh refinement levels.  

Figure~\ref{fig:bm2_cond_number} presents the condition number of the $\Mat{S'}$ for various choices of stabilization constant.  The minimum is achieved close to the recommended value $\tau^*$ that was inferred from the single macroelement analysis.  Furthermore, the removal of near-singular modes from the system matrix has a dramatic impact on iterative solver performance.  Figure~\ref{fig:bm02_krylov} presents the number of Krylov iterations to convergence needed for different mesh refinements. For low values of the stabilization constant, the near-singular modes cause a dramatic degradation in the linear solver performance.
  
\subsubsection{Modified Undrained Barry-Mercer} \label{sec:bm_3}

This last variant of the Barry Mercer's problem tests the efficacy of the stabilization when even more severe pressure oscillations are present. The setup is based on \cite{haga2012}, where the only difference with the previous setup is that the source point is switched from rate-controlled to pressure-controlled. The applied pressure $p_s$ varies according to
\begin{equation}
p_{s}(t) = p_\text{max} \sin(t),
\end{equation}
with $p_\text{max}=1$. The other simulation parameters follow \cite{haga2012} and are listed in Table~\ref{tab:parameters}. An exact solution is not available in this case. However, at early times ($t=0.01$ s) we can infer that the whole domain should have zero pressure except for the cell where the pressure is enforced. Figure \ref{fig:bm3_pressure} illustrates the solution obtained from the unstabilized and the stabilized schemes for a domain discretized with cell size $\Delta h = 1/16$. As expected, the stabilization considering $\tau = \tau^*$ eliminates the wild oscillations. Table~\ref{tab:bm03_evalue} reports the extremal eigenvalues and condition number of the $\Mat{S'}$ for solutions with and without stabilization. One observes that as the mesh is refined the minimal eigenvalue and the condition number converge to an asymptotic value different than zero and infinity, respectively, only when using the stabilized scheme. Figure~\ref{fig:bm03_cond_number} presents the condition number of the $\Mat{S'}$ matrix as a function of the stabilization constant. Once again, the minimimal condition number is attained when $\tau = \tau^*$. Table~ Figure~\ref{fig:bm03_krylov} shows the resulting improvement in Krylov convergence.

 \subsection{Multiphase Example}
 \label{sec:multiphase}

 \begin{figure}[t]
 \centering
 \includegraphics[width=0.5\textwidth]{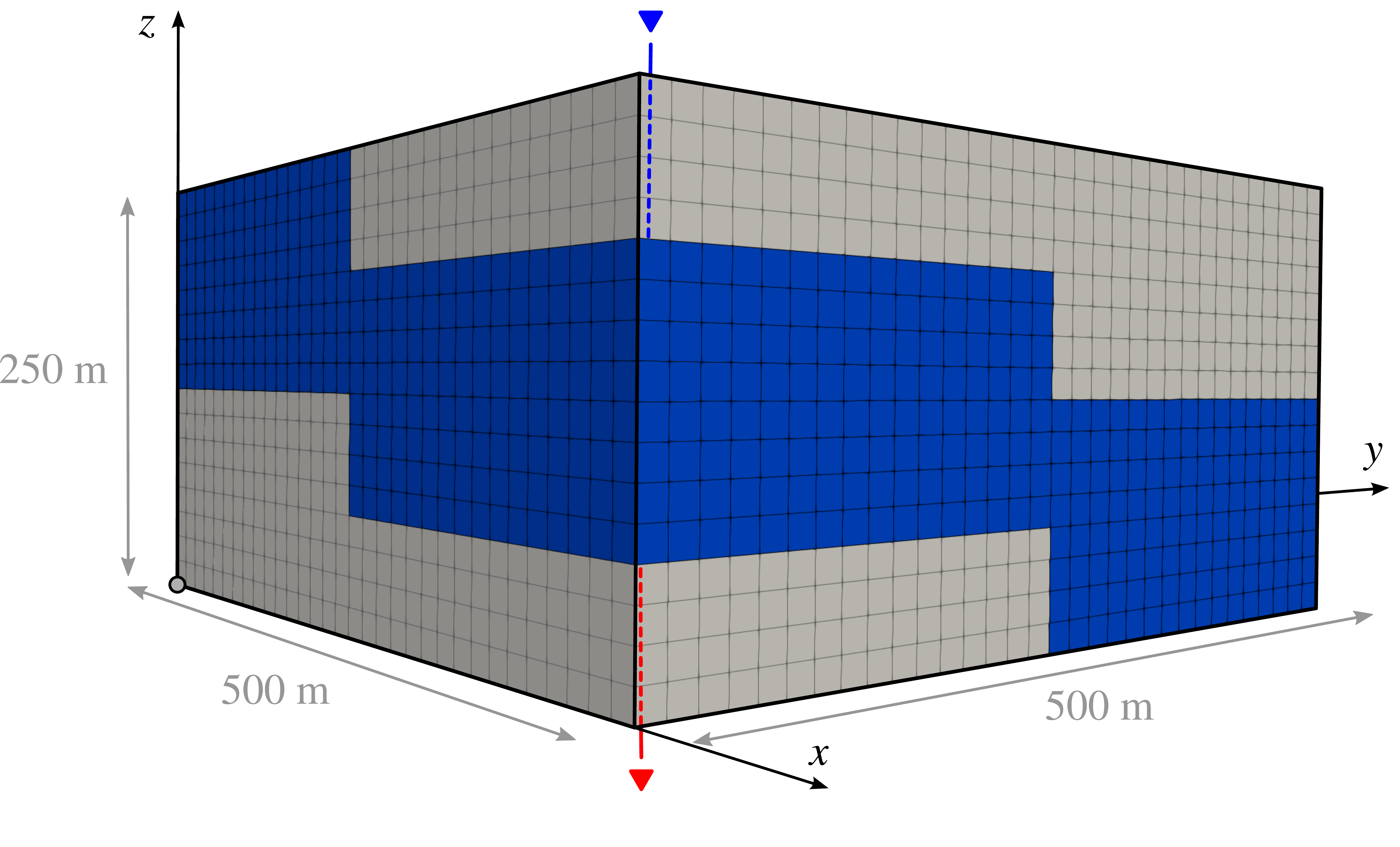}
 \caption{Problem geometry for the multiphase poromechanics example.  Grey region corresponds to a high-permeability channel, which spirals in a ``staircase'' fashion from the upper injection well (blue triangle) to the lower production well (red triangle). Blue region is a low-permeability zone.}
 \label{fig:multiphase_domain}
 \end{figure}
 
 \begin{table}[t]
\caption{Simulation parameters used in the multiphase example.}
\label{table:multiphase}
\centering
\small
\begin{tabular}{llr}
\toprule
Parameter & Units & Value  \\
\midrule 
\it Porosity: \\
$\qquad$High-perm zone & -- & 0.20  \\
$\qquad$Low-perm zone & -- & 0.05   \\[5pt]
\it Permeability: \\
$\qquad$High-perm zone & mD & 1000\\
$\qquad$Low-perm zone & mD & 1  \\[5pt]
\it Relative perm: \\
$\qquad$Residual wetting sat. & -- & 0.2  \\
$\qquad$Residual non-wetting sat. & -- & 0.2 \\[5pt]

\it Wetting Fluid:\\
$\qquad$Reference density &kg/m$^3$& 1035  \\
$\qquad$Bulk modulus &MPa& $\infty$  \\
$\qquad$Viscosity & cP & 0.3   \\[5pt]


\it Non-wetting Fluid:\\
$\qquad$Reference density &kg/m$^3$& 863  \\
$\qquad$Bulk modulus &MPa & $\infty$  \\
$\qquad$Viscosity & cP & 3.0   \\[5pt]

\it Rock: \\
$\qquad$Young's modulus &MPa& 5000  \\
$\qquad$Biot coefficient & -- & 1   \\
$\qquad$Grain density & kg/m$^3$ & 2650\\[5pt]

\it Well control:\\
$\qquad$Injection $\Delta$BHP &MPa& $5$ \\
$\qquad$Production $\Delta$BHP &MPa& $-5$ \\
$\qquad$Ramp time & day & 1 \\
$\qquad$Well radius & m & 0.1524 \\
$\qquad$Skin factor & -- & 0 \\[5pt]
\it Time-stepping: \\
$\qquad$Initial $\Delta t$ & day & 0.0001  \\
$\qquad$Maximum $\Delta t$ & day & 1 \\
$\qquad$End time & day & 100 \\[5pt]

\it Solver tolerances: \\
$\qquad$Newton & -- & 10$^{-6}$  \\
$\qquad$Krylov  & -- & 10$^{-10}$  \\
\bottomrule
\end{tabular}
\end{table}

 \begin{figure}[t]
 \begin{center}
 \begin{subfigure}[b]{\textwidth} \centering
 \includegraphics[width=0.49\textwidth]{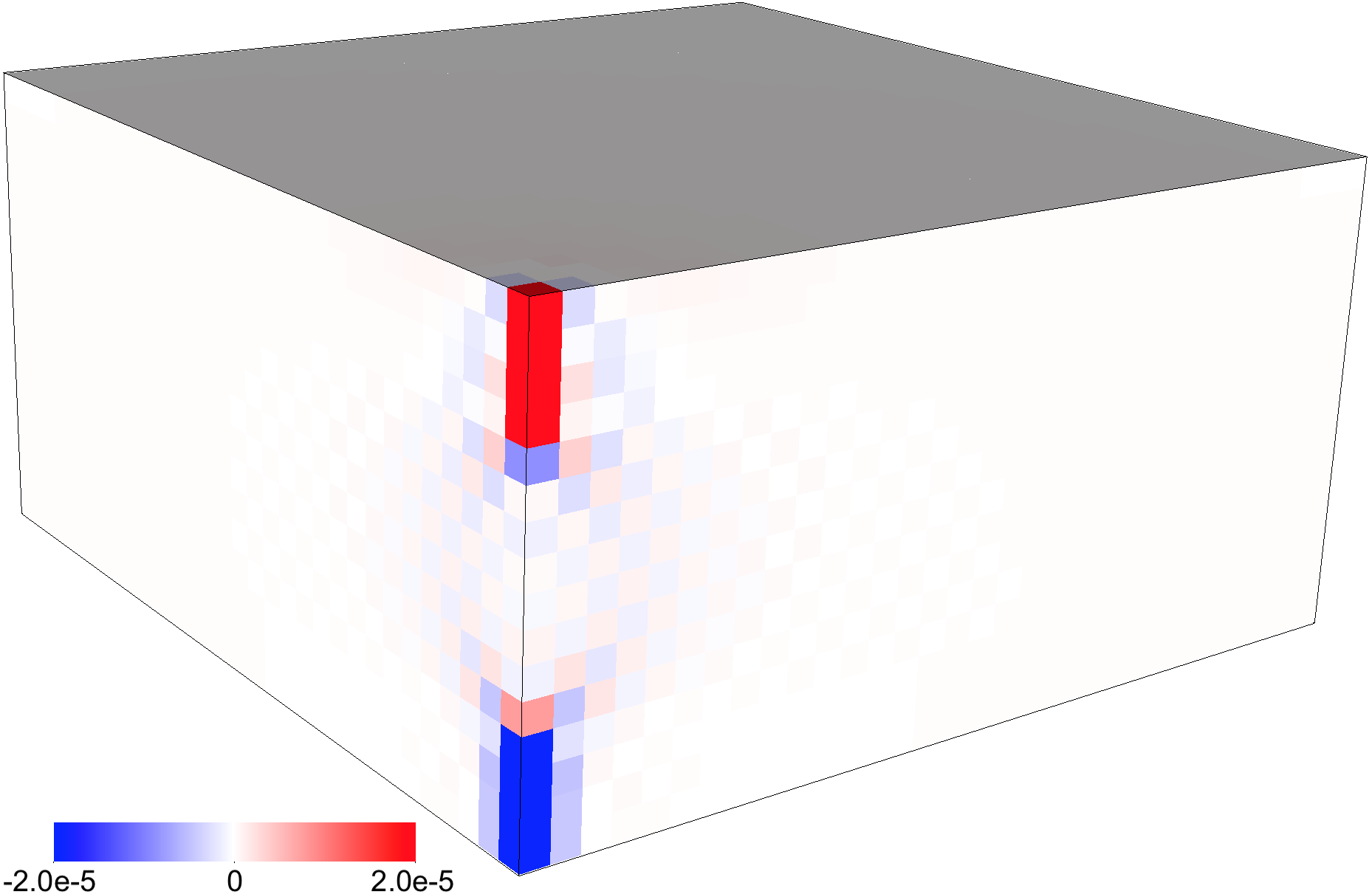}
 \includegraphics[width=0.49\textwidth]{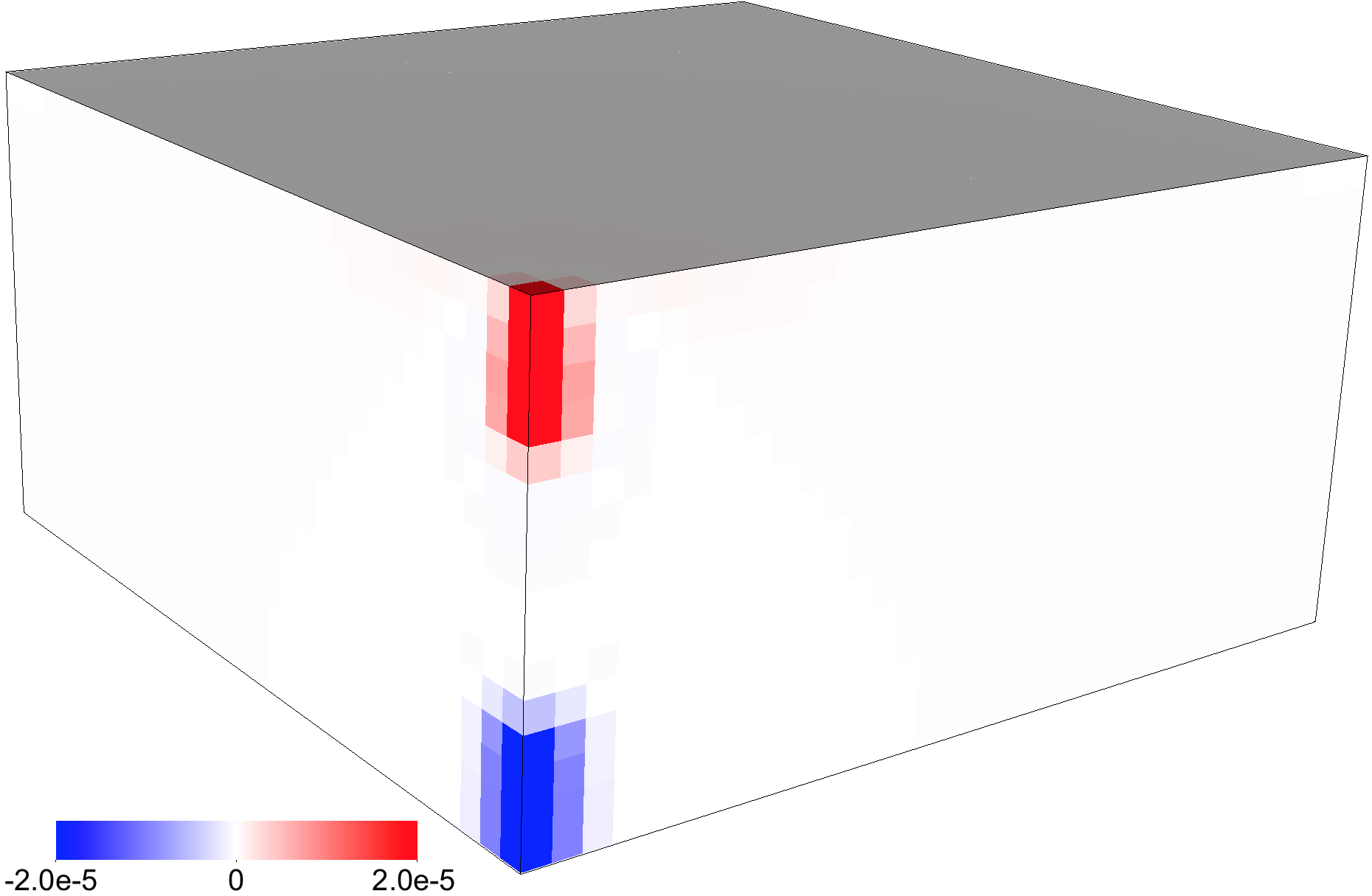}
 \caption{Colorbar is pressure in MPa at $t=10^{-4}$ day for the unstabilized (left) and stabilized (right) formulations. Note that the colorbar is truncated to visually accentuate oscillations.}
 \end{subfigure} \\
 \vspace{0.5in}
  \begin{subfigure}[b]{\textwidth} \centering
 \includegraphics[width=0.49\textwidth]{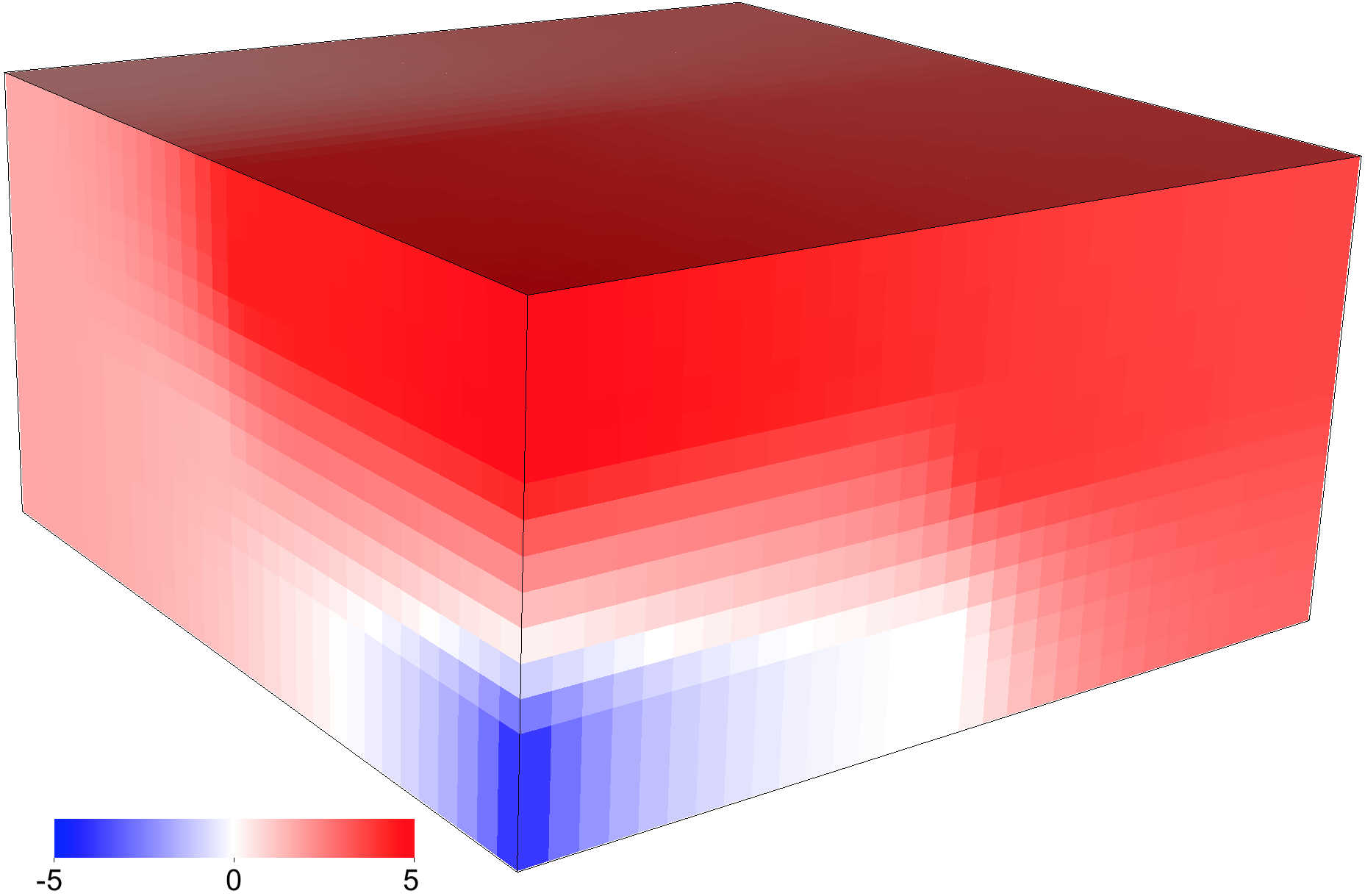}
 \includegraphics[width=0.49\textwidth]{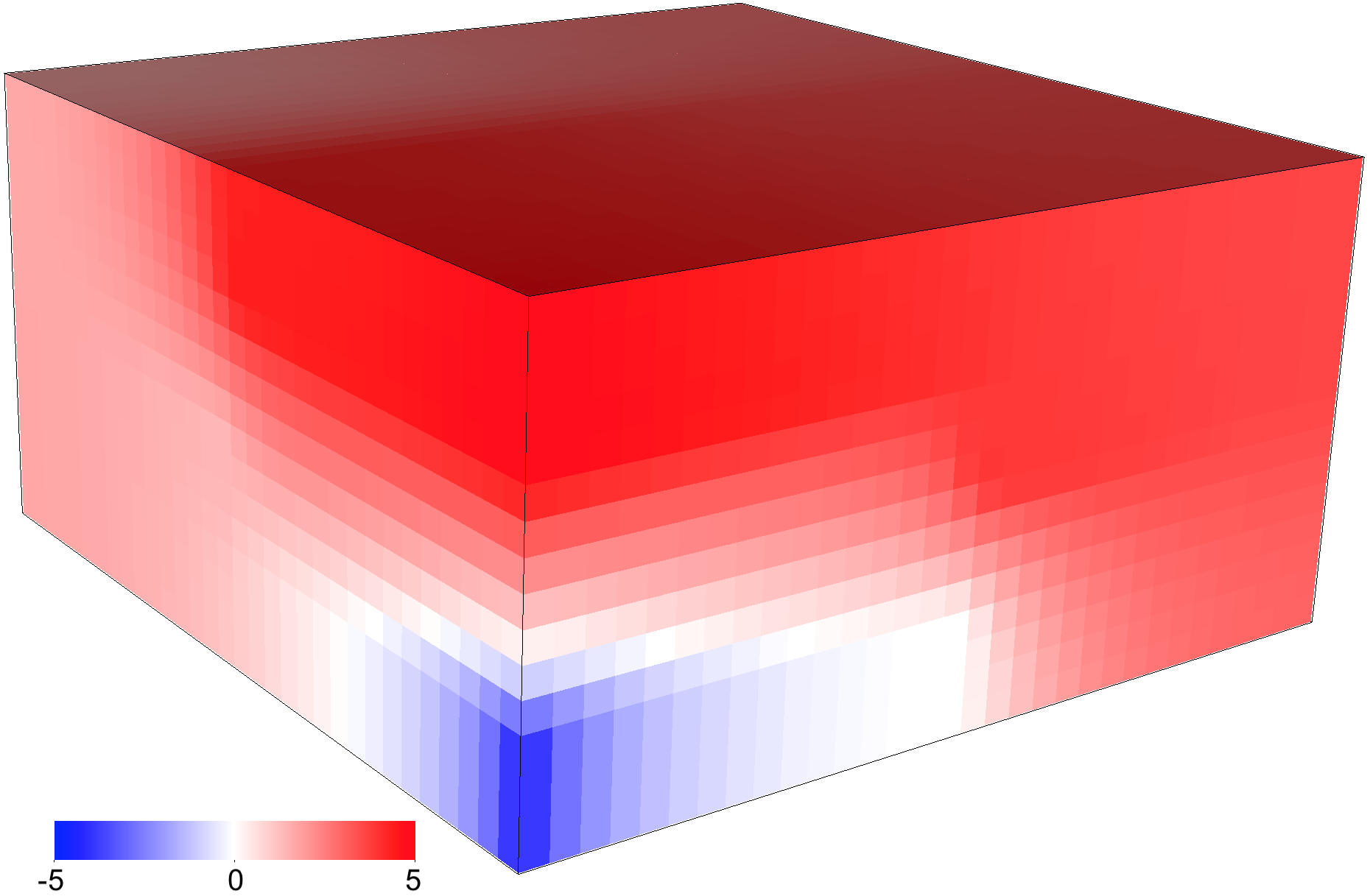}
  \caption{Colorbar is pressure in MPa at $t=100$ day for the unstabilized (left) and stabilized (right) formulations.}
 \end{subfigure}\\
 \vspace{0.5in}
  \begin{subfigure}[b]{\textwidth} \centering
 \includegraphics[width=0.49\textwidth]{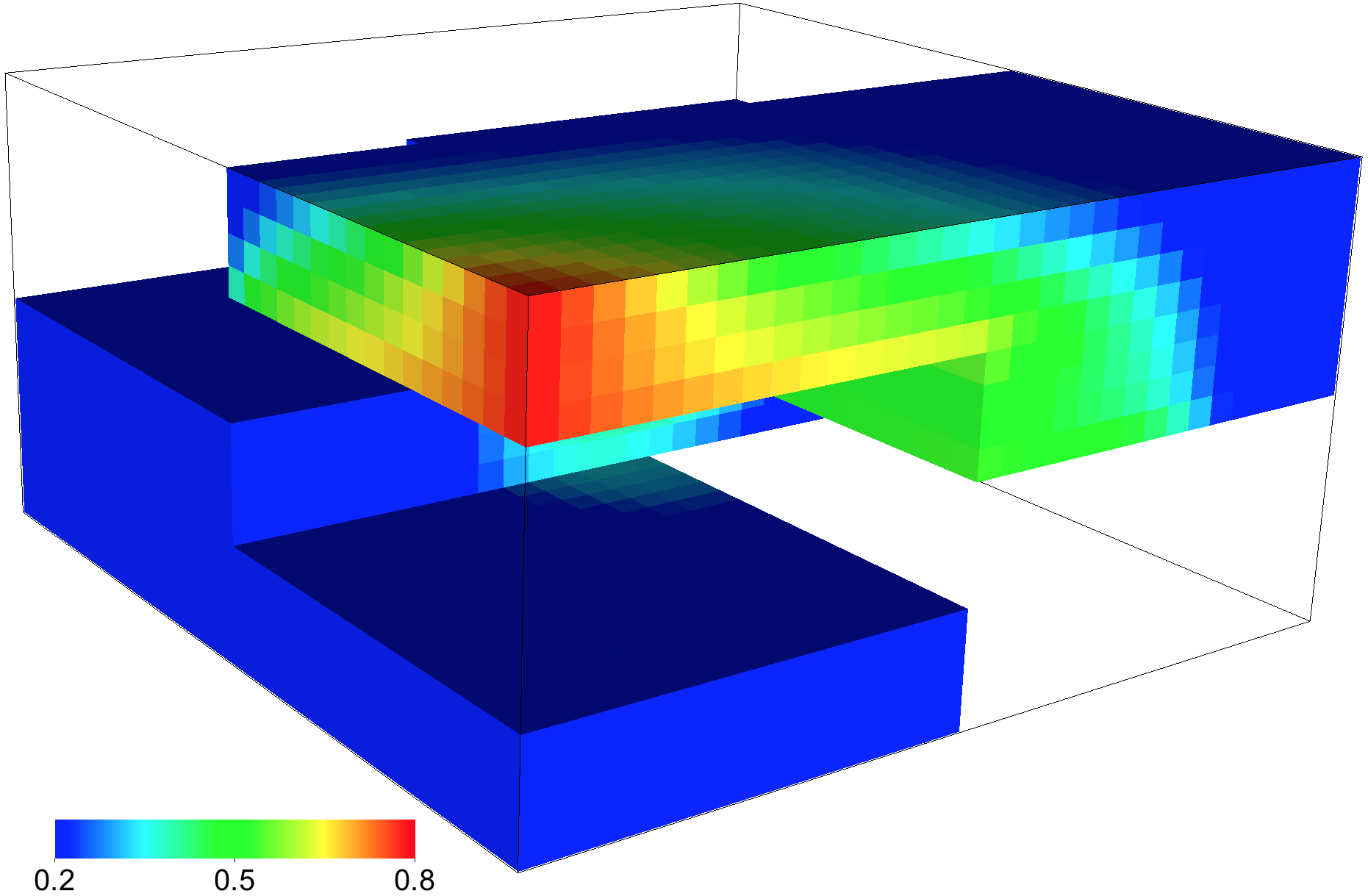}
 \includegraphics[width=0.49\textwidth]{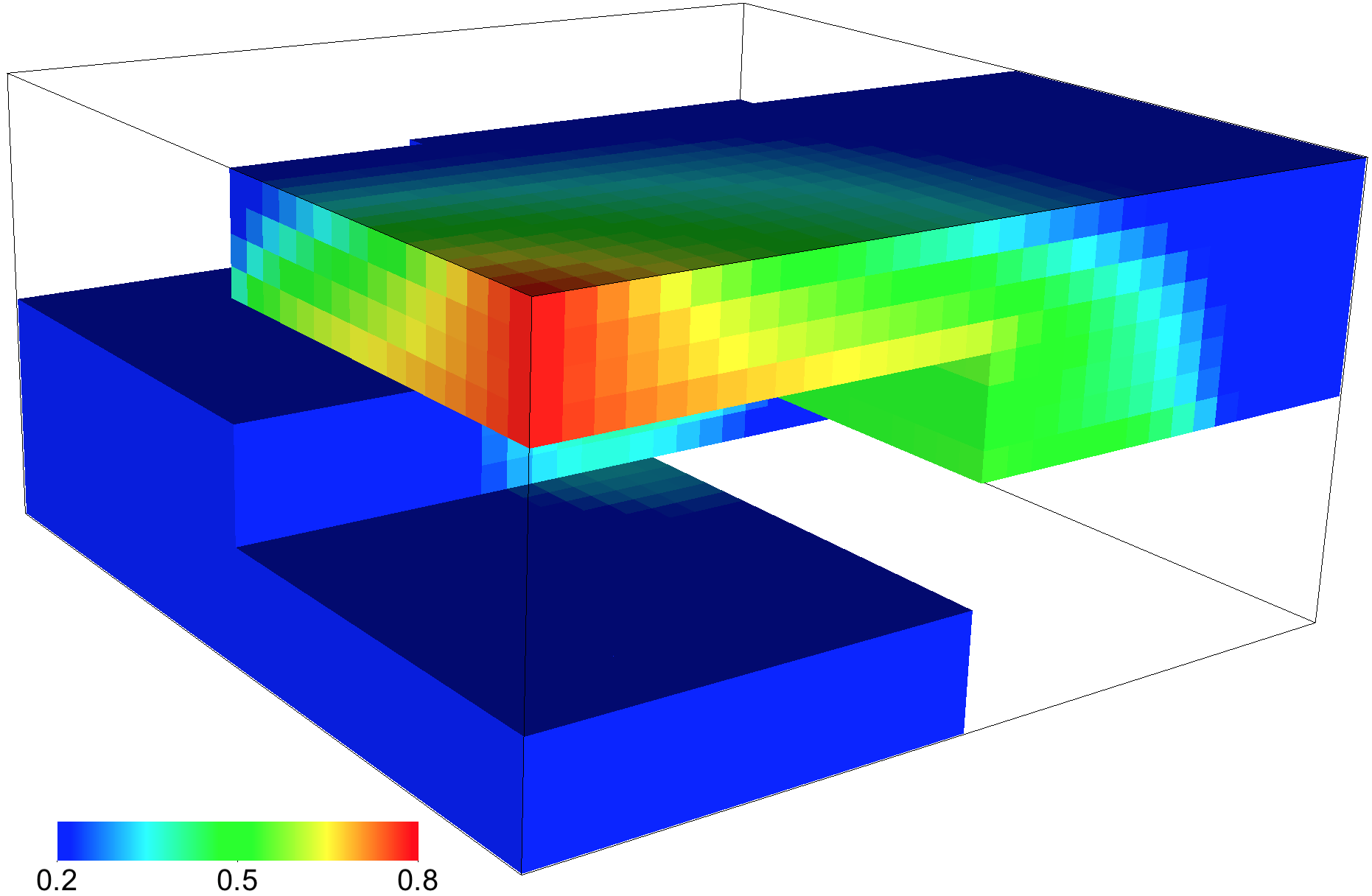}
  \caption{Colorbar is saturation in decimals at $t=100$ day for the unstabilized (left) and stabilized (right) formulations. }
 \end{subfigure}
 \caption{Comparison of unstabilized and stabilized formulations for the multiphase example.  The stabilization suppresses checkerboarding at early simulation times (a), but does not otherwise compromise solution quality at late times (b-c).}'
  \label{fig:multiphase_results}
 \end{center}
 \end{figure}
 
 \begin{figure}[t]
 \centering
 \includegraphics[width=0.45\textwidth]{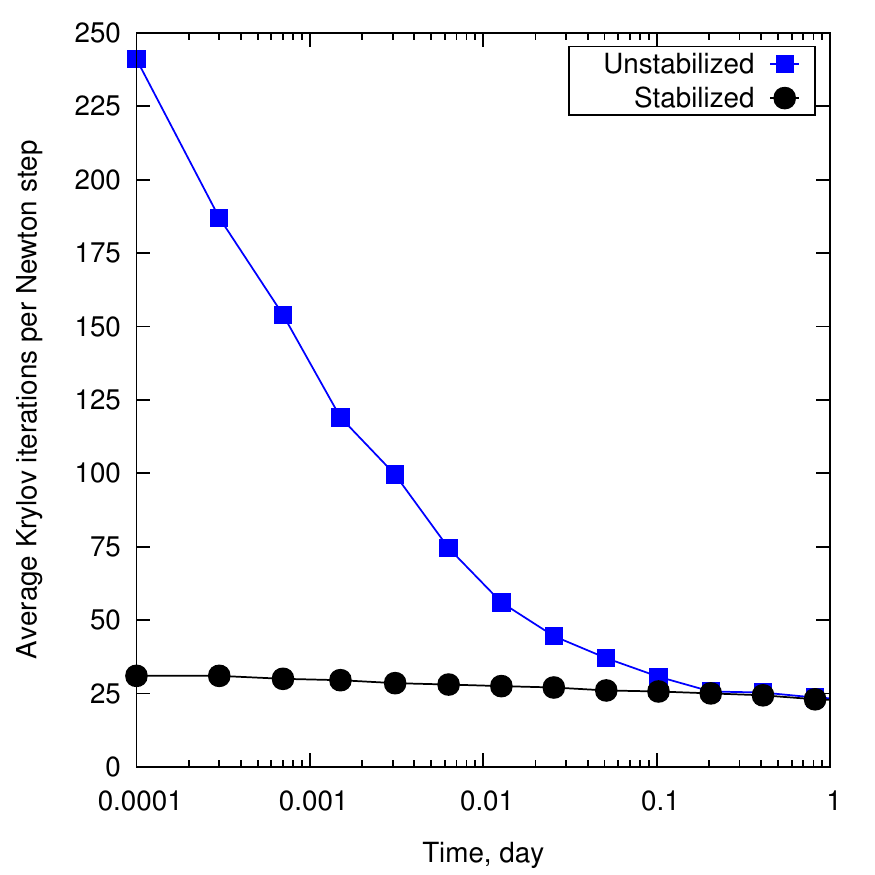}
 \caption{Comparison of average Krylov iterations per Newton step for the stabilized and unstabilized formulation applied to the multiphase poromechanics example.  Without stabilization, the iterative solver convergence degrades at early simulation times when small time steps are employed.}
  \label{fig:multiphase_krylov}
 \end{figure}

 Lastly, we consider a multiphase poromechanics example.   The test problem is based on the staircase benchmark problem originally presented in \cite{white2018}.  Figure~\ref{fig:multiphase_domain} provides an illustration of the problem geometry.  The domain contains two regions, a high-permeability channel and a low-permeability host rock.  The high-permeability channel winds its way in a spiral, staircase fashion from an upper injection well to a lower production well.  This spiral geometry is obviously artificial, but it introduces very strong coupling between the displacement, pressure, and saturation fields.  Water is injected at the upper well, while both fluids may be produced from the lower well.  These wells have a bottom-hole pressure (BHP) control.  The injector (producer) is ramped up to $5$ MPa ($-5$ MPa) overpressure over one day, and then held at a constant pressure.  All problem parameters are given in Table~\ref{table:multiphase}.  We have set the fluids to be incompressible to accentuate any instabilities in the formulation.   Specifics regarding the well control, relative permeability model, and mechanical boundary conditions may be found in \cite{white2018}. At the beginning of the simulation, the initial time step is $\Delta t = 0.0001$ day.  This time step is then doubled every step until a maximum time step of $\Delta t=1$ day is reached.  The whole simulation is run for 100 days.  Note that small time steps are the most problematic from a stability point of view.

Figure~\ref{fig:multiphase_results} presents pressure and saturation snapshots for this simulation, using both an unstabilized and stabilized formulation (with $\tau = \tau^*$).  At the first time step ($t=0.0001$ day) checkerboard oscillations are apparent in the pressure field for the unstabilized formulations.  Note that we have truncated the colorbar, cutting off the peak pressures, in order to accentuate these oscillations visually.  The stabilization successfully suppresses this checkerboarding.   At the end of the simulation ($t=100$ day), we see that the stabilized and unstabilized formulations produce essentially identical results.  The addition of the artificial flux terms does not compromise overall solution quality, with the saturation field being advected the same distance along the high-perm channel in both cases.

It is interesting to examine the linear solver behavior at early times in the simulation (Figure~\ref{fig:multiphase_krylov}).  At each Newton step of the nonlinear solver, a preconditioned GMRES iteration is used to solve the Jacobian system.  The preconditioner we use is the multistage preconditioner described in \cite{white2018}, which in general provides excellent convergence behavior for this class of problem.   In the first day of simulation time, however, we see a substantial degradation in solver performance using the unstabilized formulation.  This is a direct results of the presence of near-singular modes, to which Krylov-based solvers are extremely sensitive.  With the addition of stabilization, however, this problem is completely removed and GMRES once again exhibits excellent convergence.  Later in the simulation, as $\Delta t$ grows, the physical fluxes between elements grow and even the unstabilized formulation becomes intrinsically stable.  As a result, we see the solver convergence behavior merge at later times for the two formulations.

\section{Conclusion}
\label{sec:conclusions}
In this work, we have presented a stabilized formulation for $\mathbb{Q}_1-\mathbb{P}_0$ discretizations of single- and multiphase poromechanics. The stabilization is achieved by adding artificial flux terms to faces interior to macroelements.  We have also identified an appropriate value for the stabilization parameter based on an eigenvalue analysis of an incompressible macroelement patch test.  The stabilization is easy to implement in existing codes, and does not change the underlying sparsity pattern of the finite volume stencil.  While exact mass conservation on individual elements is sacrificed, exact mass conservation on macroelements is retained.  We have demonstrated, through a number of single and multiphase examples, that the method is effective in practice.  It can suppress spurious oscillations and also prevent unwanted degradation in iterative solver convergence in the presence of near-singular modes.  The latter is a critical issue for large-scale simulations of geosystems.

While the discussion here has been limited to $\mathbb{Q}_1-\mathbb{P}_0$ discretizations, a similar approach can likely be used for other unstable interpolation pairs involving piecewise constant pressure approximations--e.g. on more general hexahedral or tetrahedral meshes.

\section*{Acknowledgements}
Funding for JTC and JAW was provided by Total S.A. through the FC-MAELSTROM Project. JTC also acknowledges financial support provided by the Brazilian National Council for Scientific and Technological Development (CNPq) and the the John A. Blume Earthquake Engineering Center. RIB was supported by the U.S. Department of Energy, Office of Science, Office of Basic Energy Sciences, Geosciences Research Program, under Award Number DE-FG02-03ER15454. The authors wish to thank Nicola Castelletto for helpful discussions.  Portions of this work were performed under the auspices of the U.S. Department of Energy by Lawrence Livermore National Laboratory under Contract DE-AC52-07-NA27344. 

\clearpage
\newpage

\clearpage
\newpage

\bibliography{references}



\end{document}